\documentclass{article}

\usepackage[leqno,intlimits]{amsmath}
\usepackage{amssymb}
\usepackage{amsthm}
\usepackage{hyperref}
\usepackage{appendix}

\newcommand*{\be}{\begin{equation}}
\newcommand*{\ee}{\end{equation}}
\newcommand*{\ba}{\begin{aligned}}
\newcommand*{\ea}{\end{aligned}}
\newcommand*{\hop}{\bigskip\noindent}
\newcommand*{\Zb}{\mathbb Z}
\newcommand*{\un}[1]{\underline{#1}}
\newcommand*{\om}{\omega}
\newcommand*{\ze}{\zeta}
\newcommand*{\de}{\delta}
\newcommand*{\vp}{\varphi}
\newcommand*{\Rb}{\mathbb R}
\newcommand*{\e}[1]{\text{\rm e}^{#1}}
\newcommand*{\te}{\theta}
\newcommand*{\Pv}{{\bf P}}
\newcommand*{\Vv}{{\text{\bf Var}}}
\newcommand*{\Cov}{{\text{\bf Cov}}}
\newcommand*{\Ev}{{\bf E}}
\newcommand*{\vr}{\varrho}
\newcommand*{\Bc}{\mathcal B}
\newcommand*{\Fc}{\mathcal F}
\newcommand*{\Hc}{\mathcal H}
\newcommand*{\wih}{\widehat}
\newcommand*{\uomm}{\underline{\omega}^-}
\newcommand*{\uomp}{\underline\omega}
\newcommand*{\omp}{\omega}
\newcommand*{\omm}{\omega^-}
\newcommand*{\uetap}{\underline{\eta}^+}
\newcommand*{\la}{\lambda}
\newcommand*{\ga}{\gamma}
\newcommand*{\al}{\alpha}
\newcommand*{\fl}[1]{\lfloor{#1}\rfloor}
\newcommand*{\ce}[1]{\lceil{#1}\rceil}
\newcommand*{\wt}{\widetilde}
\newcommand*{\Mc}{\mathcal M}
\newcommand*{\mush}{{\un\mu}^{\text{shock}\,\vr}}
\newcommand*{\rate}{f}
\newcommand*{\abs}[1]{\lvert#1\rvert}
\newcommand*{\lrt}{(\ell,\,\err,\,\te)}
\newcommand*{\err}{\mathfrak r}
\newcommand*{\hhgt}{\un{h}}
\newcommand*{\lr}{[\ell,\,\err]}
\newcommand*{\ghgt}{\un{g}}
\newcommand*{\Nconf}{\un{N}}
\newcommand*{\abt}{(a,\,b,\,\te)}

\newcommand*{\di}{\,\text{\rm d}}

\newcommand*{\ind}{\mathbf{1}}
\newcommand*{\pcount}{I}

\newtheorem{tm}{Theorem}[section]
\newtheorem{ass}[tm]{Assumption}
\newtheorem{cor}[tm]{Corollary}
\newtheorem{lm}[tm]{Lemma}

\numberwithin{equation}{section}

\begin{document}
\title{Fluctuation bounds in the exponential bricklayers process}
\author{
 \begin{tabular}{@{}c@{\ \,}}
  M\'arton Bal\'azs\thanks{Budapest University of Technology and Economics. Part of this work was done while M.\ B.\ was affiliated with the MTA-BME Stochastics Research Group},\\
  {\small\tt balazs@math.bme.hu}
 \end{tabular}
 \begin{tabular}{@{}c@{\ \,}}
  J\'ulia Komj\'athy\thanks{Budapest University of Technology and Economics},\\
  {\small\tt komyju@math.bme.hu}
 \end{tabular}
 \begin{tabular}{@{}c@{\ \,}}
  Timo Sepp\"al\"ainen\thanks{University of Wisconsin-Madison\newline
  M.\ Bal\'azs and J.\ Komj\'athy were partially supported by the Hungarian Scientific Research Fund (OTKA) grants K60708, TS49835, F67729, and the Morgan Stanley Mathematical Modeling Center.
  M.\ Bal\'azs was also supported by the Bolyai Scholarship of the Hungarian Academy of Sciences and by the Hungarian Scientific Research Fund (OTKA) grant K100473.
  T.\ Sepp\"al\"ainen was partially supported by National Science Foundation grants DMS-0701091 and DMS-10-03651, and by the Wisconsin Alumni Research Foundation.}\\
 {\small\tt seppalai@math.wisc.edu}
 \end{tabular}
}

\maketitle
\begin{abstract}
This paper is the continuation of our earlier paper \cite{unipq3}, where we proved \(t^{1/3}\)-order of current fluctuations across the characteristics in a class of one dimensional interacting systems with one conserved quantity. We also claimed two models with concave hydrodynamic flux which satisfied the assumptions which made our proof work. In the present note we show that the \emph{totally asymmetric exponential bricklayers process} also satisfies these assumptions. Hence this is the first example with \emph{convex} hydrodynamics of a model with \(t^{1/3}\)-order current fluctuations across the characteristics. As such, it further supports the idea of universality regarding this scaling.
\end{abstract}


\hop

\section{Introduction}

It is conjectured that particle current through the characteristics of one dimensional stochastic interacting systems with one conserved quantity and concave or convex hydrodynamics show \(t^{1/3}\)-order fluctuations and Tracy-Widom type limit distributions in this order. Our earlier paper \cite{unipq3} provides a robust argument that proves this order of the fluctuations. We refer to that paper for the general framework and other results of the field. Very briefly, \cite{unipq3} works if one proves the following properties of a model (see the exact formulation therein):
\begin{enumerate}
\item a strict domination of a second class particle of a denser system on one of a sparser system,\label{it:yz}
\item a non-strict, but tight, domination of a second class particle on a system of second class particles that are defined between the system in question and another system with a different density,\label{it:qx}
\item strictly concave or convex, in the second derivative sense, hydrodynamic flux function of the hyperbolic conservation law obtained by the Eulerian limiting procedure,\label{it:hcvx}
\item a tail bound of a second class particle in a(n essentially) stationary process.\label{it:qtb}
\end{enumerate}
Properties \ref{it:yz} and \ref{it:qx} form what we call the \emph{microscopic concavity or convexity property}. Arguments in \cite{unipq3} are worked out for the concave setting, but everything works word-for-word in the convex case. 

We try to very briefly indicate here how the argument \cite{unipq3} proceeds from these assumptions. We consider two coupled processes with different densities, and a density of second class particles between them. We tag one of these second class particles, let us denote its position by $X$. We also consider the position $Q$ of a single second class particle that evolves on one of the two coupled processes. Microscopic concavity allows us a (non strict, but tight) comparison between the positions $X$ and $Q$. Briefly, the argument now goes as follows. A deviation event of $Q$ implies, via the comparison, that of $X$. Next we relate the deviation of $X$ to a surface growth deviation by simply noticing that current of second class particles is just difference in surface growth. The probability of the surface growth deviation is controlled by the variance. An exact connection \cite{varj2nd} between the variance of the surface growth and the first centered moment of $Q$ then ``closes the loop'', and now we have a deviation bound of $Q$ in terms of its first centered moment. If we close one of the two densities to the other in the right pace as a function of time, then constants appear in a way that allows to conclude the $\text{time}^{1/3}$ scaling; this is the point where strict convexity of the hydrodynamic flux function is essential.

Two examples are claimed in \cite{unipq3}: the asymmetric simple exclusion process and a totally asymmetric zero range process with jump rates that increase with exponentially decaying slope. In this note we prove the above properties and hence the \(t^{1/3}\) scaling for yet another system, the totally asymmetric exponential bricklayers process (TAEBLP). This model was introduced in \cite{valak}, and its normal fluctuations off-characteristics were demonstrated in \cite{fluct} (in case of general convex jump rates, not only exponential).

General convex increasing rates of a totally asymmetric bricklayers process allow couplings that prove properties \ref{it:yz} and \ref{it:hcvx} above. The exponential jump rates have a strong enough convexity property that will allow us to show property \ref{it:qx}.
We do this by repeating an argument somewhat similar to the one applied to the concave zero range process in \cite{unipq3}. The idea resembles much to the concave case, but we include this convex case in full details (rather than listing all the differences from previous work) due to the complexity of the method.

Finally, property \ref{it:qtb} is highly nontrivial when the jump rates have unbounded increments. We use a coupling based on property \ref{it:yz} and a recent result \cite{rwshscp} that asserts that a second class particle of the exponential bricklayers process performs a simple (drifted) random walk under appropriate shock initial conditions. It is worth noting that exponential jump rates were also of fundamental importance in \cite{rwshscp}, this technical point being the main reason for considering this particular family of jump rates in this note. Indeed, this is the only point that prevents us from proving the result for e.g.~the totally asymmetric zero range process with \emph{convex} exponential jump rates.

We emphasize at this point that we only consider nearest neighbor models. We believe that, as far as the hydrodynamic flux is strictly convex in the second derivative sense, the $t^{1/3}$ scaling should hold for a wide class of non nearest neighbor dynamics as well. However, as intricate couplings and orderings of second class particles play a crucial role in the methods, we do not see an easy way to deal with the non-nearest neighbor case.

Let us also have a few comments on explicit product invariant stationary distributions. In \cite{unipq3} we explicitly use them, as they make the arguments easier. The crucial points of the method are properties 2 (microscopic convexity) and 4 (tail bound of the second class particle) above. These depend on the details of the models, and the few known examples for which they could be proved indeed have product stationary distributions. Therefore we have not investigated how the arguments in \cite{unipq3} and the present note could be generalized to the case of other types of stationary distributions. We again believe that once microscopic convexity and the tail bound were proved, the remainder of the argument could be generalized and the scaling would remain valid for many models with non product stationary distributions as well.

The case of exponential jump rates was constructed in \cite{exists}. The results of the note \cite{varj2nd} are used by \cite{unipq3}. Those require strong construction results which are not provided by \cite{exists} and therefore, to our knowledge, are not available. To close that gap, we reproduce the results of \cite{varj2nd} here for the TAEBLP.

The organization of this paper is the following: we repeat the introduction of the model, the fluctuation results and conclusions of \cite{unipq3}, and the definition of the microscopic convexity property in Section \ref{sc:repeat}. We construct the four process coupling and prove the microscopic convexity property in Section \ref{sc:mcbl}. Finally we show how to use \cite{rwshscp} to prove property \ref{it:qtb} in Section \ref{sc:qtb}. The result of the note \cite{varj2nd} is reproduced in the Appendix.

\section{The model, properties and results}\label{sc:repeat}

\subsection{The model}

The model we discuss is the totally asymmetric exponential bricklayers process (TAEBLP) introduced in \cite{valak}, and also treated in \cite{sokvalak} and \cite{rwshscp}. The model is a member of the class in \cite{unipq3}, here is a brief definition. The process describes the growth of a surface which we imagine as the top of a wall formed by columns of bricks over the interval $(i,\,i+1)$ for each pair of neighboring sites $i$ and $i+1$ of $\Zb$. The height $h_i$ of this column is integer-valued. The components of a configuration $\un{\om}\in\Omega$ are  the negative discrete gradients of the heights: $\om_i=h_{i-1}-h_i\,\in \Zb$. The configuration space is therefore
\[
\Omega:\,=\left\{\un\om=(\om_i)_{i\in\Zb}\ :\ \om_i\in\Zb\right\}=\Zb^{\Zb}.
\]
At times it will be convenient to have notation for the increment configuration
$\un\de_i\in\Omega$ with exactly one nonzero entry equal to $1$:
\be
(\un\de_i)_j=\left\{\ba
&1,&&\text{for }i=j,\\
&0,&&\text{for }i\ne j.
\ea\right.\label{eq:dedef}
\ee

Bricklayers processes are characterized by a function \(f\,:\,\Zb\to\Rb^+\). We only consider the totally asymmetric nearest neighbor case here, in which only deposition of bricks in the following way are allowed:
\be
\left.
\ba
\left(\om_i,\,\om_{i+1}\right)&\longrightarrow\left(\om_i-1,\,\om_{i+1}+1\right)\\
h_i&\longrightarrow h_i+1
\ea
\right\}
\text{with rate}\ f(\om_i)+f(-\om_{i+1}).\label{eq:add}
\ee
Conditionally on the present state, these moves happen independently at all sites $i$. We can summarize this information in the  formal infinitesimal generator $L$ of the process $\un\om(\cdot)$:
\[
(L\vp)(\un\om)=\sum_{i\in\Zb}[f(\om_i)+f(-\om_{i+1})]\cdot\left[\vp(\un\om^{(i,\,i+1)} )-\vp(\un\om)\right].
\]
A jump results in the new configuration $\un\om^{(i,\,i+1)}$ defined  by
\[
\bigl(\un\om^{(i,\,i+1)}\bigr)_j=
\left\{\begin{array}{ll}
\om_j&\text{for}\ j\neq i,\,i+1,\\ 
\om_i-1&\text{for}\ j=i,\\
\om_{i+1}+1&\text{for}\ j=i+1.
\end{array}\right.
\]
$L$ acts on bounded cylinder functions  $\vp\,:\,\Omega\to\Rb$ (this means that $\vp$ depends only on finitely many  $\om_i$-values).
 The additive form of the rates gives rise to the bricklayers representation: at each site \(i\) stands a bricklayer who places a brick on the column on his left with rate \(f(-\om_i)\) and on the one on his right with rate \(f(\om_i)\).

Thus we have a Markov process $\{\un\om(t): t\in\Rb^+\}$ of an evolving increment configuration and a Markov process $\{\un h(t): t\in\Rb^+\}$ of an evolving height configuration. The initial increments $\un\om(0)$ specify the initial height $\un h(0)$ up to a vertical translation. We shall always normalize the height process so that \(h_0(0)=0\).

Attractivity of the process is essential for this paper. This is achieved by assuming that \(f\) is nondecreasing.

Finally, stationary translation-invariant product distributions for \(\un\om(\cdot)\) are ensured by \(f(z)\cdot f(1-z)=1\) for each \(z\in\Zb\).

The totally asymmetric \emph{exponential} bricklayers process (TAEBLP) is obtained by taking
\be
f(z)=\e{\beta(z-1/2)}.\label{eq:fdef}
\ee

The construction of the bricklayers process with any nondecreasing \(f\) that is bounded by an exponential function is given in \cite{exists} on a set of tempered configurations \(\wt\Omega\). This set consists of configurations with bounded asymptotic slope, the precise definition is given in \cite{exists}. As certain desired semigroup properties are not fully proved, we avoid technical difficulties in the proofs of \cite{varj2nd} by reproducing its results for the TAEBLP in the Appendix.

\subsection{The basic coupling}\label{sc:bc}

We use a particularly simple form of the basic coupling which is made possible by the bricklayer representation: it is enough to define the structure of moves as described in \cite{unipq3} for a given side (left or right) of an individual bricklayer. Here is how to do it for a given bricklayer at site \(i\). Given the present configurations \(\un\om^1,\,\un\om^2,\,\dots,\,\un\om^n\in\wt\Omega\), let $m\mapsto\ell(m)$ be a permutation that orders the \(\om_i\) values:
\[
\om^{\ell(m)}_i\le\om^{\ell(m+1)}_i,\qquad1\le m<n.
\]
For simplicity, set
\[
p(m):\,=f(\om^{\ell(m)}_i)\qquad\text{and}\qquad q(m):\,=f(-\om^{\ell(m)}_i),
\]
and the dummy variables \(p(0)=q(n+1)=0\). Recall that the function \(f\) is nondecreasing. Now the rule is that independently for each \(m=1,\dotsc, n\), at rate $p(m)-p(m-1)$,  precisely bricklayers of $\un\om^{\ell(m)}$, $\un\om^{\ell(m+1)}$, $\dots$, $\un\om^{\ell(n)}$ place a brick on their right, and bricklayers of \(\un\om^{\ell(1)},\,\un\om^{\ell(2)},\,\dots,\,\un\om^{\ell(m-1)}\) do not. Independently, at rate \(q(m)-q(m+1)\), precisely bricklayers of $\un\om^{\ell(1)}$, $\un\om^{\ell(2)}$, $\dots$, $\un\om^{\ell(m)}$ place a brick on their left, and bricklayers of \(\un\om^{\ell(m+1)},\,\un\om^{\ell(m+2)},\,\dots,\,\un\om^{\ell(n)}\) do not. Given the configurations \(\un\om^1,\,\un\om^2,\,\dots,\,\un\om^n\in\wt\Omega\), bricklayers at different sites perform the above steps independently.

The combined effect of these joint rates creates the correct marginal rates, that is, the bricklayer of $\un\om^{\ell(m)}$ executes the move \eqref{eq:add} with rate $p(m)=f(\om^{\ell(m)}_i)$, and the same move on column \(h_{i-1}\) with rate $q(m)=f(-\om^{\ell(m)}_i)$. 
  
Notice also that, due to monotonicity of \(f\), a jump of \(\un\om^a\) without \(\un\om^b\) on the column \([i,\,i+1]\) by the bricklayers at site \(i\) can only occur if \(f(\om^b_i)<f(\om^a_i)\) which implies \(\om^a_i>\om^b_i\). Also, a jump of \(\un\om^a\) without \(\un\om^b\) on the column \([i-1,\,i]\) by the bricklayers at site \(i\) can only occur if \(f(-\om^b_i)<f(-\om^a_i)\) which implies \(\om^a_i<\om^b_i\). The result of any of these steps then cannot increase the number of discrepancies  between the two processes, hence the name \emph{attractivity} for monotonicity of \(f\). In particular, a sitewise ordering $\om^a_i\le \om^b_i$\quad\(\forall i\in\Zb\) is preserved by the basic coupling. 

The differences between two processes are called \emph{second class particles}. Their number is nonincreasing. In particular, if, for processes \(\un\om^a\) and \(\un\om^b\) we have \(\om^a_i\ge\om^b_i\) for each \(i\in\Zb\), then the second class particles between them are conserved. A special case that is of key importance to us is the situation where only one second class particle is present between two processes.

\subsection{Hydrodynamics and some exact identities}

From now on, we restrict our attention to the TAEBLP. Recall the jump rates \eqref{eq:add}. As described in \cite{unipq3}, the process has product translation-invariant stationary distribution with marginals \(\mu^\te\)
\be
 \mu^\te(z)=\frac{\e{\te z}}{f(z)!}\cdot\frac1{Z(\te)}\label{eq:mudef}
\ee
that turn out to be of discrete Gaussian type, see \cite{valak} for the explicit formula. \(\Pv^\te,\ \Ev^\te,\ \Vv^\te,\ \Cov^\te\) will refer to laws of a process evolving in this stationary distribution. The density \(\vr(\te):\,=\Ev^\te(\om)\in\Rb\) is a strictly increasing function of the parameter \(\te\in\Rb\), and can take on any real value by the Appendix of \cite{unipq3}. \(\mu^\vr,\ \Pv^\vr,\ \Ev^\vr,\ \Vv^\vr,\ \Cov^\vr\) will refer to laws of a density \(\vr\) stationary process.

The hydrodynamic flux is
\[
\Hc(\vr)=\Ev^\vr [f(\om)+f(-\om)]=\e{\te(\vr)}+\e{-\te(\vr)}.
\]
As \(f\) \eqref{eq:fdef} is convex and nonlinear, the Appendix of \cite{unipq3} applies and yields a convex hydrodynamic flux with
\be
\Hc''(\vr)>0.\label{eq:hconv}
\ee
(This convexity property is quite natural, as $\Hc$ is the expected value of the convex rate function $f$. It has a nontrivial proof \cite{unipq3} that uses total positivity or, a proof \cite{convex} using correlation inequalities.)

As in \cite{unipq3}, we introduce
\be
\wih\mu^\vr(y):\,=\frac{1}{\Vv^\vr(\om_0)}\sum_{z=y+1}^{\infty}(z-\vr)\mu^\vr(z),\qquad y\in\Zb.\label{eq:muhat}
\ee
The Appendix of \cite{unipq3} applies to show that both \(\mu^\vr\) and \(\wih\mu^\vr\) are stochastically monotone in \(\vr\). Denote by \(\Ev\) the expectation w.r.t.\ the evolution of a pair \((\uomm(\cdot),\,\uomp(\cdot))\) started with initial data (recall \eqref{eq:dedef})
\be
\uomm(0)=\uomp(0)-\un\de_0\sim\Bigl(\bigotimes_{i\ne0}\mu^\vr\Bigr)\otimes\wih\mu^\vr,\label{eq:umhdef}
\ee
and evolving under the basic coupling. This pair will always have a single second class particle whose position is denoted by \(Q(t)\).  In other words, $\uomm(t)=\uomp(t)-\un\de_{Q(t)}$. We reprove Corollaries 2.4 and 2.5 of \cite{varj2nd} in the Appendix that state that for any \(i\in\Zb\) and \(t\ge0\),
\begin{align}
\Vv^\vr(h_i(t))&=\Vv^\vr(\om)\cdot\Ev|Q(t)-i|\label{eq:varcovar}
\intertext{and}
\Ev(Q(t))&=V^\vr\cdot t,\notag
\end{align}
where \(V^\vr=\Hc'(\vr)\) is the characteristic speed.  Note in particular that in \eqref{eq:varcovar} the variances are taken in a stationary process, while the expectation of $Q(t)$ is taken in the coupling with initial distribution \eqref{eq:umhdef}. 

\subsection{Results}

We repeat the results of \cite{unipq3}, valid now for the TAEBLP. 

\begin{tm}
Fix any density $\vr\in\Rb$, let the TAEBLP processes $(\uomm(t),\,\uomp(t))$ e\-vol\-ve in basic coupling with
initial distribution {\rm\eqref{eq:umhdef}} and let $Q(t)$ be the position of the second class
particle between $\uomm(t)$ and $\uomp(t)$.
Then there is a constant \(C_1=C_1(\vr)\in(0,\,\infty)\) such that for all \(1\le m<3\),
\[
\frac1{C_1}<\liminf_{t\to\infty}\frac{\Ev|Q(t)-V^\vr t|^m}{t^{2m/3}}\le\limsup_{t\to\infty}\frac{\Ev|Q(t)-V^\vr t|^m}{t^{2m/3}}<\frac{C_1}{3-m},
\]
where $\Ev$ on the right hand-side refers to the expectation in a coupled pair started from almost-equlibrium \eqref{eq:umhdef}.
\end{tm}
Superdiffusivity of the second class particle is best seen with the choice \(m=2\): the variance of its position is of order \(t^{4/3}\).
Next some corollaries. 
Notation \(\fl{X}\) stands for the lower integer part of \(X\).
\begin{cor}[Current variance]
There is a constant \(C_1=C_1(\vr)>0\), such that
\[
\frac1{C_1}<\liminf_{t\to\infty}\frac{\Vv^\vr(h_{\fl{V^\vr t}}(t))}{t^{2/3}}\le\limsup_{t\to\infty}\frac{\Vv^\vr(h_{\fl{V^\vr t}}(t))}{t^{2/3}}<C_1.
\]
\end{cor}
\noindent
\begin{cor}[Weak Law of Large Numbers for the second class particle]\label{cr:lln}
In a density-\(\vr\) stationary process,
\[
\frac{Q(t)}{t}\overset{\text{d}}{\to}V^\vr.
\]
\end{cor}
\begin{cor}[Dependence of current on the initial configuration]
For any \(V\in\Rb\) and \(\al>1/3\) the following limit holds in the \(L^2\) sense for a density-\(\vr\) stationary process:
\be
\lim_{t\to\infty}\frac{h_{\fl{Vt}}(t)-h_{\fl{Vt}-\fl{V^\vr t}}(0)-t(\Hc(\vr)-\vr\Hc'(\vr))}{t^\al}=0.
\label{eq:L2cor}\ee
\end{cor}
Recall that
\[
h_{\fl{Vt}-\fl{V^\vr t}}(0)=\left\{\ba
&\sum_{i=\fl{Vt}-\fl{V^\vr t}+1}^0\om_i(0),&&\text{if }V<V^\vr,\\
&\qquad\quad0,&&\text{if }V=V^\vr,\\
&-\sum_{i=1}^{\fl{Vt}-\fl{V^\vr t}}\om_i(0),&&\text{if }V>V^\vr
\ea\right.
\]
only depends on a finite segment of the initial configuration. Limit
\eqref{eq:L2cor} shows that on the 
diffusive time scale $t^{1/2}$  only fluctuations from the initial distribution are visible: these
fluctuations are  translated rigidly at the characteristic speed \(V^\vr\). 
\begin{cor}[Central Limit Theorem for the current]
For any \(V\in\Rb\) in a density-\(\vr\) stationary process
\[
\lim_{t\to\infty}\frac{\Vv^\vr(h_{\fl{Vt}}(t))}{t}=\Vv^\vr(\om)\cdot|V^\vr-V|=\,:D,
\]
and the Central Limit Theorem also holds: the centered and normalized  height
 \(\wt h_{\fl{Vt}}(t)/\sqrt{t\cdot D}\) converges in distribution to a standard normal.
\label{cr:CLT}\end{cor}
\noindent

For convex rate zero range and bricklayers processes Corollaries \ref{cr:lln} and \ref{cr:CLT} were proved  by M.\ Bal\'azs \cite{fluct}.  

The way we obtain the Theorem and the corollaries is simply proving that the assumptions formulated in \cite{unipq3} hold. For the sake of completeness, we repeat these, translated to our convex case.

\subsection{Microscopic convexity}

We start with the definition of \emph{microscopic convexity}. This is just translated from the microscopic concavity property of \cite{unipq3}, where more detailed explanations and comments can be found. Convexity \eqref{eq:hconv} implies that the characteristic speed $V^\vr=\Hc'(\vr)$ 
is a nondecreasing function of the density $\vr$:
\be
\la<\vr\;\Longrightarrow V^\la\le V^\vr.\label{eq:Vorder}
\ee
The microscopic counterpart 
of  a characteristic is the motion of a second class particle.  Our key assumption
that we term microscopic convexity is that the ordering \eqref{eq:Vorder} 
can also be realized at the particle level as an ordering between two second class 
particles introduced into two processes at densities $\la$ and $\vr$.

Let \(\la<\vr\) be two densities.
Define \({\wih\mu}^\vr+1\) as the measure that gives weight
\({\wih\mu}^\vr(z-1)\) to an integer $z\in\Zb$ (recall \eqref{eq:muhat}).
By the stochastic domination  \({\wih\mu}^\la\le{\wih\mu}^\vr\), we can
let \({\wih\mu}^{\la,\vr}\) be a coupling measure with
marginals \({\wih\mu}^\la\) and \({\wih\mu}^\vr+1\) and with the property
\[
{\wih\mu}^{\la,\vr}\{(y,\,z)\,:\,y<z\}=1.
\]
Let also
\(\mu^{\la,\vr}\) be a coupling measure of site-marginals \(\mu^\la\)
and \(\mu^\vr\) of the invariant distributions,  with
\be
\mu^{\la,\vr}\{(y,\,z)\,:\,y\le z\}=1.\label{eq:cmu}
\ee
Note the distinction that under ${\wih\mu}^{\la,\vr}$ the second coordinate
is strictly above the first. 

To have notation for inhomogeneous product measures on $\Zb^\Zb$, let 
\(\un\la=(\la_i)_{i\in\Zb}\) and \(\un\vr=(\vr_i)_{i\in\Zb}\) denote sequences
of density values, with  \(\la_i\) and \(\vr_i\) assigned to site \(i\).  
The product distribution with marginals \({\wih\mu}^{\la_0,\vr_0}\)
at the origin and \(\mu^{\la_i,\vr_i}\) at other sites is denoted by
\[
{\un{\wih\mu}}^{\un\la,\un\vr}:\,=\Bigl(\bigotimes_{i\ne0}\mu^{\la_i,\vr_i}\Bigr)\otimes{\wih\mu}^{\la_0,\vr_0}.
\]
Measure \({\un{\wih\mu}}^{\un\la,\un\vr}\) gives probability one to the event
\[
\{(\un\eta(0),\,\uomp(0))\,:\,\eta_0(0)<\omp_0(0)\text{, and }\eta_i(0)
\le\omp_i(0)\text{ for }0\ne i\in\Zb\}.
\]
The initial configuration \((\un\eta(0),\,\uomp(0))\) will always be assumed a member
of this set, and the pair process \((\un\eta(t),\,\uomp(t))\) evolves in basic coupling. 
Notice that \({\un{\wih\mu}}^{\un\la,\un\vr}\) is in general \emph{not} stationary for this joint evolution.

The discrepancies between these two processes are called the
\emph{ \(\omp-\eta\) (second class) particles}. The number of such particles at site
 $i$ at time $t$ is $\omp_i(t)-\eta_i(t)$.   In the basic coupling   the
 \(\omp-\eta\) particles are
conserved, in the sense that none are created or annihilated.
 We label the  \(\omp-\eta\) particles with
 integers, and let $X_m(t)$ denote the position of particle $m$ at time $t$. 
 The  initial 
labeling is chosen to   satisfy
\[
\dotsm\le X_{-1}(0)\le X_0(0)=0<X_{1}(0)\le\dotsm.
\]
We can specify that $X_0(0)=0$ because under  \({\un{\wih\mu}}^{\un\la,\un\vr}\)
there is an \(\omp-\eta\) particle at site $0$ with probability $1$. 
During the evolution we keep the positions $X_i(t)$  of the  \(\omp-\eta\) particles ordered.
To achieve this we stipulate that 
\be\begin{array}{l}
\text{whenever an  \(\omp-\eta\) particle jumps from
a site, }\\
\text{if the jump is to the right  the highest label moves,}\\
\text{and if the jump is to the left  the  lowest label moves.}\end{array}\label{eq:Xrule}\ee

Here is the precise form of microscopic convexity for this paper.  The assumption 
states that a certain joint construction of processes (that is, a coupling) can be
performed for a range of densities in a neighborhood of a fixed density $\vr$. 
Recall \eqref{eq:dedef} for the definition of the configuration \(\un\de\).

\begin{ass}\label{def:microc}
Given a density $\vr\in\Rb$, there exists $\ga_0>0$ such that the  following holds. 
For any  $\un\la$ and $\un\vr$ such that  \(\vr-\ga_0\le \la_i\le\vr_i\le\vr+\ga_0\) for all \(i\in\Zb\), 
 a joint process  \((\un\eta(t),\,\uomp(t),\,y(t),\,z(t))_{t\ge 0}\) can be constructed
with the following properties.  
\begin{itemize} 
\item Initially \((\un\eta(0),\,\uomp(0))\) is \({\un{\wih\mu}}^{\un\la,\un\vr}\)-distributed
and the joint process  \((\un\eta(\cdot),\,\uomp(\cdot))\) evolves in basic coupling.
\item Processes $y(\cdot)$ and  $z(\cdot)$ are integer-valued. Initially \(y(0)=z(0)=0\).
  With probability one
\be
y(t)\ge z(t)\text{ for all }t\ge0.\label{eq:yz}
\ee
\item Define the processes 
\be
\uomm(t):\,=\uomp(t)-\un\de_{X_{y(t)}(t)}\quad\text{and}\quad\uetap(t):\,=\un\eta(t)+\un\de_{X_{z(t)}(t)}.\label{eq:indef}
\ee
Then both pairs $(\un\eta,\,\uetap)$ and  $(\uomm,\,\uomp)$  evolve marginally in basic
coupling. 
\item  For each $\ga\in(0,\,\ga_0)$ and large enough $t\ge 0$ 
 there exists a probability  distribution \(\nu^{\vr,\ga}(t)\) 
on \(\Zb^+\)  satisfying the tail bound
\be
\nu^{\vr,\ga}(t)\{y\,:\,y\ge y_0\}\le Ct^{\kappa-1}\ga^{2\kappa-3}y_0^{-\kappa}\label{eq:nuc}
\ee
for some fixed constants  \(3/2\le\kappa<3\) and $C<\infty$, and 
 such that  if   \(\vr-\ga\le \la_i\le\vr_i\le\vr+\ga\) for all \(i\in\Zb\), 
 then we have the stochastic bounds 
 \be y(t)\overset{\text{d}}{\ge}-\nu^{\vr,\ga}(t)\quad\text{and}\quad 
z(t)\overset{\text{d}}{\le}\nu^{\vr,\ga}(t).  \label{eq:yzbds}\ee
\end{itemize}
\end{ass}

Let us clarify some of the details in this assumption. 

Equation \eqref{eq:indef} says that \(Q^\eta(t):\,=X_{z(t)}(t)\) is  the single
second class particle between \(\un\eta\) and \(\uetap\),
 while \(Q(t):\,=X_{y(t)}(t)\) is  the one between \(\uomm\) and \(\uomp\). 
The first three bullets   say that it is possible to construct jointly four processes 
 \((\un\eta,\,\uetap,\,\uomm,\,\uomp)\) with the  specified initial conditions and so that 
 each pair 
  \((\un\eta,\,\uomp)\), $(\un\eta,\,\uetap)$ and $(\uomm,\,\uomp)$ has
   the desired marginal distribution,
  and most importantly so that 
\[
Q^\eta(t)=X_{z(t)}(t)\le X_{y(t)}(t)=Q(t).
\]
 This is a consequence of \eqref{eq:yz}  because the 
  \(\omp-\eta\) particles $X_i(t)$  stay ordered.

The tail bound \eqref{eq:nuc} is formulated in this somewhat complicated fashion because
this appears  to be the weakest form our present proof allows.  For the TAEBLP
$\nu^{\vr,\ga}(t)$ will actually be a fixed geometric distribution.

The assumptions made imply \(\un\eta(t)\le\uomp(t)\) a.s.,
and by   \eqref{eq:indef}  
\[
\un\eta(t)\le\uetap(t)\le\uomp(t)\quad\text{and}\quad\un\eta(t)\le\uomm(t)\le\uomp(t)\quad\text{a.s.}
\]
In our actual construction for the TAEBLP it turns out that while the triples \((\un\eta,\,\uetap,\,\uomp)\) and \((\un\eta,\,\uomm,\,\uomp)\) evolve also in
 basic coupling, the full joint evolution \((\un\eta,\,\uetap,\,\uomm,\,\uomp)\) does not.

As already explained, the microscopic convexity idea is contained in 
inequality \eqref{eq:yz}. There is also a sense in which 
the tail bounds \eqref{eq:yzbds} relate to convexity of the flux. 
Consider the situation
  \(\la_i\equiv\la<\vr\equiv\vr_i\).  We would expect the \(\omp-\eta\) particle \(X_0(\cdot)\)
to have average and long-term velocity
\[
R(\la,\,\vr)=\frac{\Hc(\vr)-\Hc(\la)}{\vr-\la},
\]
the Rankine-Hugoniot or shock speed. By convexity 
 \(\Hc'(\vr)=V^\vr\ge R(\la,\,\vr)\ge V^\la=\Hc'(\la)\). A strict microscopic counterpart 
 would be \(y(t)\ge0\ge z(t)\). But this condition is overly restrictive.  The 
 distributional bounds \eqref{eq:yzbds} are the natural relaxations we use.

Section \ref{sc:mcbl} contains the proof of Assumption \ref{def:microc} for the TAEBLP. The proof of \eqref{eq:yzbds} makes use of the particular exponential form \eqref{eq:fdef} of the rates. Unfortunately, we do not have an argument for more general convex rates at the moment.

There is one more assumption in \cite{unipq3} needed to state the main result. Constants \(C_{\centerdot},\,\al_{\centerdot}\) will not depend on time, but might depend on the density parameter \(\vr\), and their values can change from line to line.
\begin{ass}\label{ass:main}
Let $(\uomm,\,\uomp)$ be a pair of processes in basic coupling,
started from distribution \eqref{eq:umhdef}, with second class particle $Q(t)$. 
Then there exist constants \(0<\al_0,\,C<\infty\) such that 
\[
\Pv\{|Q(t)|>K\}\le
C\cdot\frac{t^2}{K^3}
\]
whenever \(K>\al_0 t\) and $t$ is large enough.
\end{ass}
Such an assumption is natural and easy to prove if the jump rates have bounded increments. Since \(f\) \eqref{eq:fdef} does not, this statement for the TAEBLP is nontrivial. We prove it in Section \ref{sc:qtb} for the TAEBLP.

\section{Proof of microscopic convexity}\label{sc:mcbl}

In this section we verify that Assumption \ref{def:microc} can be satisfied.
  The task is to construct   the processes \(y(t)\) and \(z(t)\) with the requisite properties.
 First let the processes \((\un\eta(\cdot),\,\uomp(\cdot))\) 
 evolve in the basic coupling so that \(\eta_i(t)\le\omp_i(t)\) for all $i\in\Zb$ and $t\ge 0$.
  We consider as a background process this pair with the labeled and ordered  $\omp-\eta$ 
 second class particles \(\dotsm\le X_{-2}(t)\le X_{-1}(t)\le X_0(t)\le X_1(t)\le X_2(t)\le\dotsm\).
  
  At each time \(t\ge0\) this background induces a partition \(\{\Mc_i(t)\}\) of the label
  space \(\Zb\)
  into intervals  indexed by sites \(i\in\Zb\), with partition intervals given by 
\[
\Mc_i(t):\,=\{m\,:\,X_m(t)=i\}.
\]
(For simplicity we assumed infinitely many second class particles in both directions, but no problem arises in case we only have finitely many of them.)   \(\Mc_i(t)\) contains the labels of 
the second class particles that reside at site \(i\) at time \(t\), and can be empty. 
The labels of the second class particles that are at the same site as the one labeled \(m\) 
form the set \(\Mc_{X_m(t)}(t)=\,:\{a^m(t),\,a^m(t)+1,\,\dots,\,b^m(t)\}\). 
The processes $a^m(t)$ and $b^m(t)$ are always well-defined and satisfy
$a^m(t)\le m\le b^m(t)$. Notice that
\be
|\Mc_{X_m(t)}(t)|=b^m(t)-a^m(t)+1=\omp_{X_m(t)}(t)-\eta_{X_m(t)}(t).\label{eq:meo}
\ee

Let us clarify these notions by discussing the ways in which \(a^m(t)\) and \(b^m(t)\) can change.
\begin{itemize}
\item A second class particle jumps from site \(X_m(t-)-1\) to site \(X_m(t-)\). Then this one necessarily has label \(a^m(t-)-1\), and it becomes the lowest labeled one at site \(X_m(t-)=X_m(t)\) after the jump. Hence \(a^m(t)=a^m(t-)-1\).
\item A second class particle jumps from site \(X_m(t-)+1\) to site \(X_m(t-)\). Then this one necessarily has label \(b^m(t-)+1\), and it becomes the highest labeled one at site \(X_m(t-)=X_m(t)\) after the jump. Hence \(b^m(t)=b^m(t-)+1\).
\item A second class particle, different from \(X_m\), jumps from site \(X_m(t-)\) to site \(X_m(t-)+1\). Then this one is necessarily labeled \(b^m(t-)\), and it leaves the site \(X_m(t-)\), hence \(b^m(t)=b^m(t-)-1\).
\item A second class particle, different from \(X_m\), jumps from site \(X_m(t-)\) to site \(X_m(t-)-1\). Then this one is necessarily labeled \(a^m(t-)\), and it leaves the site \(X_m(t-)\), hence \(a^m(t)=a^m(t-)+1\).
\item The second class particle \(X_m\) is the highest labeled on its site, that is, \(m=b^m(t-)\), and it jumps to site \(X_m(t-)+1\). Then this particle becomes the lowest labeled in the set \(\Mc_{X_m(t-)+1}=\Mc_{X_m(t)}\), hence \(a^m(t)=m\). In this case \(b^m(t)\) can be computed from \eqref{eq:meo}, the number of second class particles at the site of \(X_m\) after the jump.
\item The second class particle \(X_m\) is the lowest labeled on its site, that is, \(m=a^m(t-)\), and it jumps to site \(X_m(t-)-1\). Then this particle becomes the highest labeled in the set \(\Mc_{X_m(t-)-1}=\Mc_{X_m(t)}\), hence \(b^m(t)=m\). In this case \(a^m(t)\) can be computed from \eqref{eq:meo}, the number of second class particles at the site of \(X_m\) after the jump.
\end{itemize}
We fix initially \(y(0)=z(0)=0\).  The evolution of $(y,\,z)$ is superimposed on
the background evolution $(\un\eta,\,\uomp,\,\{X_m\})$   following 
the
general rule
below: 
Immediately after every move of the background process that involves the site where
$y$ resides before this move,   $y$ picks a new value from the labels on the site where it 
resides after the move.   
Thus $y$ itself jumps only within partition
intervals $\Mc_i$. But $y$  joins a  new partition interval   whenever it is the
highest $X$-label on its site and its 
  ``carrier''
particle $X_y$ is  forced to move to the next site on the right, or it is the
lowest $X$-label on its site and its 
  ``carrier''
particle $X_y$ is  forced to move to the next site on the left.

These are the situations when
$y(t-)=b^{y(t-)}(t-)$ and at time  $t$ an $\omp-\eta$ move from this site to the right
happens, or $y(t-)=a^{y(t-)}(t-)$ and at time  $t$ an $\omp-\eta$ move from this site to the left
happens. 
(Recall that the choice of $X$-particle to move is determined by rule \eqref{eq:Xrule}.)  All this works for $z$ in exactly the same way.  

Next we specify the probabilities that $y$ and $z$ use to refresh their values. Recall \eqref{eq:fdef}. To simplify notation, we abbreviate, given integers \(\eta<\om\),
\begin{align}
p(\eta,\,\om)&=\frac{f(\omp)-f(\omp-1)}{f(\omp)-f(\eta)}=\frac{f(-\eta)-f(-\eta-1)}{f(-\eta)-f(-\omp)}=\frac{\e{\beta(\om-\eta)}-\e{\beta(\om-\eta-1)}}{\e{\beta(\om-\eta)}-1}\label{eq:pdef}\\
\intertext{and}
q(\eta,\,\om)&=\frac{f(-\omp+1)-f(-\omp)}{f(-\eta)-f(-\omp)}=\frac{f(\eta+1)-f(\eta)}{f(\omp)-f(\eta)}=\frac{\e{\beta}-1}{\e{\beta(\om-\eta)}-1}.\label{eq:qdef}
\end{align}
Notice that both \(p(\eta,\,\om)\) and \(q(\eta,\,\om)\) only depend on \(\om-\eta\). Therefore, with a little abuse of notation, we write \(p(\om-\eta):\,=p(\eta,\,\om)\), \(q(\om-\eta):\,=q(\eta,\,\om)\). Then
\[
p(1)=q(1)=1,\qquad\qquad p(d)\ge q(d),\qquad p(d)+q(d)\le1\qquad\text{for }2\le d\in\Zb.
\]
When $y$ and $z$ reside at separate sites, they refresh independently.  When they
are together in the same partition interval, they use the joint distribution 
 in the third bullet below.
\begin{itemize}
\item Whenever any change occurs in either \(\uomp\) or \(\un\eta\) at site \(X_{y(t-)}(t-)\) and, as a result of the jump,
\(a^{y(t-)}(t)\ne a^{z(t-)}(t)\), that is, \(y(t-)\) and \(z(t-)\) belong to different parts \emph{after} the jump, we abbreviate
\[
p=p(\eta_{X_{y(t-)}(t)}(t),\,\omp_{X_{y(t-)}(t)}(t))\qquad q=q(\eta_{X_{y(t-)}(t)}(t),\,\omp_{X_{y(t-)}(t)}(t))
\]
of \eqref{eq:pdef} and \eqref{eq:qdef} in the formulas below. These depend on the values of the respective processes at the site where the label \(y\) can be found right after the jump. In this case, independently of everything else,
\be
y(t):\,=\left\{\ba
&a^{y(t-)}(t),&&\text{with prob.\ }q,\\
&b^{y(t-)}(t)-1,&&\text{with prob.\ }1-p-q,\\
&b^{y(t-)}(t),&&\text{with prob.\ }p,
\ea\right.\label{eq:ych}
\ee
except for \(y(t):\,=a^{y(t-)}(t)=b^{y(t-)}(t)\) when the difference
\[
\omp_{X_{y(t-)}(t)}(t)-\eta_{X_{y(t-)}(t)}(t)
\]
is 1. Notice that the second line in \eqref{eq:ych} has probability zero when this difference is 2.
\item Whenever any change occurs in either \(\uomp\) or \(\un\eta\) at site \(X_{z(t-)}(t-)\) and, as a result of the jump,
\(a^{y(t-)}(t)\ne a^{z(t-)}(t)\), that is, \(y(t-)\) and \(z(t-)\) belong to different parts \emph{after} the jump, we abbreviate
\[
p=p(\eta_{X_{z(t-)}(t)}(t),\,\omp_{X_{z(t-)}(t)}(t))\qquad q=q(\eta_{X_{z(t-)}(t)}(t),\,\omp_{X_{z(t-)}(t)}(t))
\]
of \eqref{eq:pdef} and \eqref{eq:qdef} in the formulas below. These depend on the values of the respective processes at the site where the label \(z\) can be found right after the jump. In this case, independently of everything else,
\be
z(t):\,=\left\{\ba
&a^{z(t-)}(t),&&\text{with prob.\ }p,\\
&a^{z(t-)}(t)+1,&&\text{with prob.\ }1-p-q,\\
&b^{z(t-)}(t),&&\text{with prob.\ }q,
\ea\right.\label{eq:zch}
\ee
except for \(z(t):\,=a^{z(t-)}(t)=b^{z(t-)}(t)\) when the difference
\[
\omp_{X_{z(t-)}(t)}(t)-\eta_{X_{z(t-)}(t)}(t)
\]
is 1. Notice that the second line in \eqref{eq:zch} has probability zero when this difference is 2.
\item Whenever any change occurs in either \(\uomp\) or \(\un\eta\) at sites \(X_{y(t-)}(t-)\) or \(X_{z(t-)}(t-)\) and, as a result of the jump, \(a^{y(t-)}(t)=a^{z(t-)}(t)\), that is, \(y(t-)\) and \(z(t-)\) belong to the same part \emph{after} the jump, that is, \(X_{y(t-)}(t)=X_{z(t-)}(t)\) then we have
\[
\omp_{X_{y(t-)}(t)}(t)=\omp_{X_{z(t-)}(t)}(t)\qquad\text{and}\qquad\eta_{X_{y(t-)}(t)}(t)=\eta_{X_{z(t-)}(t)}(t),
\]
and we abbreviate
\[
p=p(\eta_{X_{y(t-)}(t)}(t),\,\omp_{X_{y(t-)}(t)}(t))\qquad q=q(\eta_{X_{y(t-)}(t)}(t),\,\omp_{X_{y(t-)}(t)}(t))
\]
of \eqref{eq:pdef} and \eqref{eq:qdef} in the formulas below. These depend on the values of the respective processes at the site where both the labels \(y\) and \(z\) can be found right after the jump. In this case, independently of everything else,
\be
\begin{pmatrix}y(t)\\z(t)\end{pmatrix}:\,=\left\{\ba
&\begin{pmatrix}a^{y(t-)}(t)\\a^{y(t-)}(t)\end{pmatrix}\text{, with prob.\ }q,\\
&\begin{pmatrix}b^{y(t-)}(t)-1\\a^{y(t-)}(t)\end{pmatrix}\text{, with prob.\ }(p-q)\land(1-p-q),\\
&\begin{pmatrix}b^{y(t-)}(t)\\a^{y(t-)}(t)\end{pmatrix}\text{, with prob.\ }[2p-1]^+,\\
&\begin{pmatrix}b^{y(t-)}(t)-1\\a^{y(t-)}(t)+1\end{pmatrix}\text{, with prob.\ }[1-2p]^+,\\
&\begin{pmatrix}b^{y(t-)}(t)\\a^{y(t-)}(t)+1\end{pmatrix}\text{, with prob.\ }(p-q)\land(1-p-q),\\
&\begin{pmatrix}b^{y(t-)}(t)\\b^{y(t-)}(t)\end{pmatrix}\text{, with prob.\ }q,
\ea\right.\label{eq:yzch}
\ee
except for \(y(t)=z(t):\,=a^{y(t-)}(t)=b^{y(t-)}(t)\) when the difference\\
\(\omp_{X_{y(t-)}(t)}(t)-\eta_{X_{y(t-)}(t)}(t)\) is 1. Notice that the second, the fourth and the fifth lines have probability zero when this difference is 2.
\end{itemize}
The above moves for \(y\) and \(z\) always occur within labels   at a given site. 
This determines whether the particle
 \(Q(t):\,=X_{y(t)}(t)\) or  \(Q^\eta(t):\,=X_{z(t)}(t)\)  is the one to jump if the next move 
 out of the site is   an $\omp-\eta$ move. 
 
We prove that the above construction has the properties required in
Assumption \ref{def:microc}. First note that the refreshing rule \eqref{eq:yzch} marginally gives the same moves and probabilities as \eqref{eq:ych} or \eqref{eq:zch} for \(y(\cdot)\) or \(z(\cdot)\), respectively.
\begin{lm}
The pair $(\uomm,\,\uomp):\,=(\uomp-\un\delta_{X_y},\,\uomp)$ obeys basic coupling, 
as does the pair $(\un\eta,\,\uetap):\,=(\un\eta,\,\un\eta+\un\delta_{X_z})$. 
\end{lm}
\begin{proof}
We write the proof for \((\uomm,\,\uomp)\). We need to show that, given the 
configuration \((\un\eta,\,\uomp,\,\{X_m\},\,y)\), the jump rates of \((\uomm,\,\uomp)\) 
are the ones prescribed in  basic coupling (Section \ref{sc:bc}) and by 
\(\eqref{eq:add}\). As mentioned in Section \ref{sc:bc}, the effect of bricklayers determine the evolution of processes. Notice first that an \(\omp-\eta\) particle can only jump away from a site \(i\) if a bricklayer of \(\omp\) or \(\eta\) moves. As the moves \eqref{eq:ych} or \eqref{eq:yzch} by themselves never result in a change of \(X_{y(\cdot)}(\cdot)\), any move of \(Q\) from a site \(i\) is a result of a bricklayer's move at site \(i\). Therefore, we see that moves initiated by bricklayers of \(\uomp\) at sites \(i\ne Q\) happen as well to \(\uomm\), as required by the basic coupling. The only point to consider  is moves by the bricklayers at site \(i=Q\). We start with them putting a brick on their right. Since the last time any change occurred at site \(i\), \(y\) chose values according to \eqref{eq:ych} or \eqref{eq:yzch}. Notice that \eqref{eq:ych} and \eqref{eq:yzch} give the same marginal probabilities for this choice. Hence 
\begin{align}
\text{\(y\) took a value less than \(b^y\) with probability}\quad1-p&=\frac{f(\omp_i-1)-f(\eta_i)}{f(\omp_i)-f(\eta_i)}\label{eq:ynjpr}
\intertext{and}
\text{\(y\) took on value \(b^y\) with probability}\quad p&=\frac{f(\omp_i)-f(\omp_i-1)}{f(\omp_i)-f(\eta_i)},\label{eq:yjpr}
\end{align}
as given in \eqref{eq:ych}. Notice that \eqref{eq:ynjpr} happens with probability zero if \(\om_i=\eta_i+1\). According to the basic coupling of \(\un\eta\) and \(\uomp\), the following right moves of bricklayers at \(i\) can occur:
\begin{itemize}
\item With rate \(f(\omp_i)-f(\eta_i)\), \(\uomp\) jumps without \(\un\eta\). The highest labeled second class particle, \(X_{b^y}\) jumps from site \(i\) to site \(i+1\).
\begin{itemize}
\item With probability \eqref{eq:yjpr} \(X_y=Q\) jumps with \(X_{b^y}\). In this case
\[
\omm_i(t-)=\omp_i(t-)-1=\omp_i(t)=\omm_i(t)
\]
since the difference \(Q\) disappears from site \(i\). Also,
\[
\omm_{i+1}(t-)=\omp_{i+1}(t-)=\omp_{i+1}(t)-1=\omm_{i+1}(t),
\]
since the difference \(Q\) appears at site \(i+1\). So in this case 
$\uomp$ undergoes a jump but    \(\uomm\) does not, and the rate is 
\[
[f(\omp_i)-f(\eta_i)]\cdot\frac{f(\omp_i)-f(\omp_i-1)}{f(\omp_i)-f(\eta_i)}=f(\omp_i)-f(\omm_i).
\]
\item With probability \eqref{eq:ynjpr} \(X_y=Q\) does not jump with \(X_{b^y}\), since it has label less than \(b^y\) (this probability is zero if \(\omp_i=\eta_i+1\)). In this case \(\uomm\) and  \(\uomp\) perform the same jump   and it occurs with rate
\[
[f(\omp_i)-f(\eta_i)]\cdot\frac{f(\omp_i-1)-f(\eta_i)}{f(\omp_i)-f(\eta_i)}=f(\omm_i)-f(\eta_i).
\]
\end{itemize}
\item With rate \(f(\eta_i)\), both bricklayers of \(\un\eta\) and \(\uomp\) at site \(i\) move. No change occurs in the \(\omp-\eta\) particles, hence no change occurs in \(Q\). This implies that the process \(\uomm\) jumps as well.
\end{itemize}
Summarizing  we see that the rate for the bricklayers of \((\uomm,\,\uomp)\) at site \(i\) to lay brick on their rights together is \(f(\om_i^-)\), and the rate for the one of $\uomp$ to move without $\uomm$ is $f(\omp_i)-f(\omm_i)$.  This is exactly what basic coupling requires.

Consider now bricklayers at site \(i=Q\) putting a brick on their left. Since the last time any change occurred at site \(i\), \(y\) chose values according to \eqref{eq:ych} or \eqref{eq:yzch}. Hence 
\begin{align}
\text{\(y\) took on value \(a^y\) with probability}\quad q&=\frac{f(-\omp_i+1)-f(-\omp_i)}{f(-\eta_i)-f(-\omp_i)}\label{eq:yla}
\intertext{and}
\text{\(y\) took a value higher than \(a^y\) with probability}\quad1-q&=\frac{f(-\eta_i)-f(-\omp_i+1)}{f(-\eta_i)-f(-\omp_i)},\label{eq:yhla}
\end{align}
as given in \eqref{eq:ych}. Notice that \eqref{eq:yhla} happens with probability zero if \(\om_i=\eta_i+1\). According to the basic coupling of \(\un\eta\) and \(\uomp\), the following left moves of bricklayers at \(i\) can occur:
\begin{itemize}
\item With rate \(f(-\eta_i)-f(-\omp_i)\), \(\un\eta\) jumps without \(\uomp\). The lowest labeled second class particle, \(X_{a^y}\) jumps from site \(i\) to site \(i-1\).
\begin{itemize}
\item With probability \eqref{eq:yla} \(X_y=Q\) jumps with \(X_{a^y}\). In this case
\[
\omm_i(t-)=\omp_i(t-)-1=\omp_i(t)-1=\omm_i(t)-1
\]
since the difference \(Q\) disappears from site \(i\). Also,
\[
\omm_{i-1}(t-)=\omp_{i-1}(t-)=\omp_{i-1}(t)=\omm_{i-1}(t)+1,
\]
since the difference \(Q\) appears at site \(i+1\). So in this case 
$\uomm$ undergoes a jump but \(\uomp\) does not, and the rate is 
\[
[f(-\eta_i)-f(-\omp_i)]\cdot\frac{f(-\omp_i+1)-f(-\omp_i)}{f(-\eta_i)-f(-\omp_i)}=f(-\omm_i)-f(-\omp_i).
\]
\item With probability \eqref{eq:yhla} \(X_y=Q\) does not jump with \(X_{a^y}\), since it has label more than \(a^y\) (this probability is zero if \(\omp_i=\eta_i+1\)). In this case none of \(\uomm\) or \(\uomp\) move; this occurs with rate
\[
[f(-\eta_i)-f(-\omp_i)]\cdot\frac{f(-\eta_i)-f(-\omp_i+1)}{f(-\eta_i)-f(\omp_i)}=f(-\eta_i)-f(-\omm_i).
\]
\end{itemize}
\item With rate \(f(-\omp_i)\), both bricklayers of \(\un\eta\) and \(\uomp\) at site \(i\) move. No change occurs in the \(\omp-\eta\) particles, hence no change occurs in \(Q\). This implies that the process \(\uomm\) jumps as well.
\end{itemize}
Summarizing  we see that the rate for the bricklayers of \((\uomm,\,\uomp)\) at site \(i\) to lay brick on their rights together is \(f(-\omp_i)\), and the rate for the one of $\uomm$ to move without $\uomp$ is $f(-\omm_i)-f(-\omp_i)$.  This is exactly what basic coupling requires. 

A very similar argument can be repeated for \((\un\eta,\,\uetap)\).
\end{proof}
\begin{lm}\label{lm:zryz}
Inequality \eqref{eq:yz} $y\ge z$ holds in the above construction.
\end{lm}
\begin{proof}
Since no jump of $y$ or $z$ moves one of them into a new partition interval, the 
only situation that can jeopardize \eqref{eq:yz} is the simultaneous refreshing 
of $y$ and $z$ in a common partition interval.  But this case is governed by 
step \eqref{eq:yzch}  which by definition ensures that  \(y\ge z\). (When \(b^{y(t-)}(t)=a^{y(t-)}(t)+1\), we have, by \eqref{eq:meo}, \(\omp_{X_{y(t-)}(t)}(t)-\eta_{X_{y(t-)}(t)}(t)=2\), and hence \(p\) of \eqref{eq:yzch} is more than 1/2. Therefore the probability of the step in line 4 of \eqref{eq:yzch} is zero.)
\end{proof}
Define the geometric distribution
\be
\nu(m):\,=\left\{\ba
&\e{-\beta m}(1-\e{-\beta}),&&m\ge0\\
&0,&&m<0.
\ea\right.\label{eq:nugeom}
\ee
 
\begin{lm}\label{lm:zrgeom}
Conditioned on the process \((\un\eta,\,\uomp)\), the bounds \(y(t)\overset{\text{d}}{\ge}-\nu\) and \(z(t)\overset{\text{d}}{\le}\nu\) hold for all $t\ge 0$.
\end{lm}
To avoid unnecessary complications with negative values, we show the proof for \(z(t)\). Notice that both the statement and the behavior of \(y(t)\) is reflected compared to \(z(t)\), hence the proof is the same for the two processes. The argument consists of three steps.
\begin{lm}\label{lm:ygeom}
Let \(Z\) be a random variable with distribution \(\nu\), and fix integers  \(a\le b\)
and  \(\eta<\omp\) so that \(\omp-\eta=b-a+1\).  Apply the following operation to $Z$: 
\begin{itemize}
\item[{\rm (i)}] if \(a\le Z\le b\), apply the probabilities from \eqref{eq:zch} with   parameters \(a,\,b,\,\eta,\,\omp\) to pick a new value for $Z$; 
\item[{\rm (ii)}] if $Z<a$ or $Z>b$ then do not change $Z$.
\end{itemize}
Then the resulting distribution \(\nu^*\) is stochastically dominated by \(\nu\).
\end{lm}
\begin{proof}
There is nothing to prove when \(b=a\), hence we assume \(b>a\) or, equivalently, \(\omp-\eta=b-a+1\ge2\). It is also clear that \(\nu^*(m)=\nu(m)\) for \(m<a\) or \(m>b\). We need to prove, in view of the distribution functions,
\[
\sum_{\ell=a}^m\nu^*(\ell)\ge\sum_{\ell=a}^m\nu(\ell)
\]
for all \(a\le m\le b\). Notice that \(\nu^*\) gives zero weight on values \(a+1<m<b\) (if any), and also that the display becomes an equality if \(m=b\). Therefore, it is enough to prove the inequality for \(m=a\):
\be
\nu^*(a)\ge\nu(a),\label{eq:nua}
\ee
and \(m=b-1\):
\be
\sum_{\ell=a}^{b-1}\nu^*(\ell)\ge\sum_{\ell=a}^{b-1}\nu(\ell)\qquad\text{that is,}\qquad\nu^*(b)\le\nu(b).\label{eq:nub}
\ee
Notice that \eqref{eq:nua} is trivially true for \(a<0\). For \(a\ge0\) we start with rewriting the left hand-side of \eqref{eq:nua} with the use of \eqref{eq:zch}, \eqref{eq:pdef}, and the abbreviation \(d=\om-\eta=b-a+1\):
\[
\ba
\nu^*(a)&=p(d)\cdot\sum_{\ell=a}^b\nu(\ell)\\
&=\frac{\e{\beta d}-\e{\beta(d-1)}}{\e{\beta d}-1}\cdot(\e{-\beta a}-\e{-\beta(b+1)})\\
&=\e{-\beta a}\cdot\frac{\e{\beta d}-\e{\beta(d-1)}}{\e{\beta d}-1}\cdot(1-\e{-\beta d})=\nu(a).
\ea
\]
As for \eqref{eq:nub}, both sides become zero if \(b<0\). For \(b\ge0\) we have
\[
\ba
\nu^*(b)&=q(d)\cdot\sum_{\ell=a}^b\nu(\ell)\\
&\le\frac{\e{\beta}-1}{\e{\beta d}-1}\cdot(\e{-\beta a}-\e{-\beta(b+1)})\\
&=\e{-\beta b}\cdot\frac{\e{\beta}-1}{\e{\beta d}-1}\cdot(\e{\beta(d-1)}-\e{-\beta})=\nu(b).
\ea
\]
\end{proof}
\begin{lm}\label{lm:yattr}
The dynamics defined by \eqref{eq:zch} is attractive.
\end{lm}
\begin{proof}
Following the same realizations of \eqref{eq:zch}, we see that two copies of \(z(\cdot)\) under a common environment can be coupled so that whenever they get to the same part \(\Mc_i\), they move together from that moment.
\end{proof}
\begin{proof}[Proof of Lemma \ref{lm:zrgeom}]
Initially \(z(0)=0\) by definition, which is clearly a distribution dominated by \(\nu\) of \eqref{eq:nugeom}. Now we argue recursively: by time \(t\) the distribution of \(z(t)\) was a.s.\ only influenced by finitely many jumps of the environment, which resulted in distributions \(\nu_1\), then \(\nu_2\), then \(\nu_3\), 
etc. Suppose \(\nu_k\overset{\text{d}}{\le}\nu\), and let \(\nu^*\) be the distribution that would result from \(\nu\) by the \(k+1^\text{st}\) jump. Then \(\nu_{k+1}\overset{\text{d}}{\le}\nu^*\) by \(\nu_k\overset{\text{d}}{\le}\nu\) and Lemma \ref{lm:yattr}, while \(\nu^*\overset{\text{d}}{\le}\nu\) by Lemma \ref{lm:ygeom}.
\end{proof}

\section{A tail bound for the second class particle}\label{sc:qtb}

In this section we prove that Assumption \ref{ass:main} holds for the TAEBLP model. The difficulty comes from the fact that jump rates of the second class particle, being the increments of the growth rates \eqref{eq:fdef}, are unbounded.
First recall the coupling measure \(\mu^{\la,\vr}\) of \eqref{eq:cmu} and notice that it gives weight one on pairs of the form \((y,\,y)\) if \(\la=\vr\). Define also \(\mu^{\text{shock}\,\vr}\) by
\[
\mu^{\text{shock}\,\vr}(y,\,z)=\left\{\ba
&\mu^\vr(y),&&\text{if }z=y+1,\\
&0,&&\text{otherwise.}
\ea\right.
\]
With these marginals we define the shock product distribution
\be
\mush:\,=\bigotimes_{i<0}\mu^{\vr+1,\vr+1}\cdot\bigotimes_{i=0}\mu^{\text{shock}\,\vr}\cdot\bigotimes_{i>0}\mu^{\vr,\vr},\label{eq:mushdf}
\ee
a measure on a pair of coupled processes with a single second class particle at the origin.
\begin{lm}
The first marginal of \(\mush\) is the product distribution
\[
\bigotimes_{i<0}\mu^{\vr+1}\cdot\bigotimes_{i\ge0}\mu^\vr,
\]
while the second marginal is
\be
\bigotimes_{i\le0}\mu^{\vr+1}\cdot\bigotimes_{i<0}\mu^\vr.\label{eq:mushsec}
\ee
\end{lm}
\begin{proof}
The first part of the statement and the second part, apart from \(i=0\), follow from the definitions. The nontrivial part is
\[
\mu^{\vr+1}(z)=\mu^\vr(z-1),\qquad z\in\Zb,
\]
valid for the second marginal at \(i=0\). This is specific to the definition \eqref{eq:mudef} of $\mu^\vr$, and of the exponential rates \eqref{eq:fdef}, and to prove it we write, with \(\te=\te(\vr)\),
\[
\mu^\vr(z-1)=\frac{f(z)}{\e{\te}}\cdot\frac{\e{\te z}}{f(z)!}\cdot\frac1{Z(\te)}=\frac{\e{(\te+\beta)z}}{f(z)!}\cdot\frac1{\e{\te+\beta/2}Z(\te)}.
\]
Summing this up for all \(z\in\Zb\) gives one on the left hand-side, hence leads to
\[
Z(\te+\beta)=\sum_{z=-\infty}^\infty\frac{\e{(\te+\beta)z}}{f(z)!}=\e{\te+\beta/2}Z(\te),
\]
which also implies
\[
\vr(\te+\beta)=\frac{\di}{\di\te}\log(Z(\te+\beta))=\vr(\te)+1.
\]
We conclude that
\[
\mu^\vr(z-1)=\frac{\e{(\te+\beta)z}}{f(z)!}\cdot\frac1{Z(\te+\beta)}=\mu^{\vr(\te+\beta)}(z)=\mu^{\vr+1}(z),
\]
which finishes the proof of the lemma.
\end{proof}

The translation of \(\mush\) is denoted by
\[
\tau_k\mush:\,=\bigotimes_{i<k}\mu^{\vr+1,\vr+1}\cdot\bigotimes_{i=k}\mu^{\text{shock}\,\vr}\cdot\bigotimes_{i>k}\mu^{\vr,\vr}.
\]
The main tool we use is Theorem 1 from \cite{rwshscp}, which we reformulate here. \(\mu S(t)\) will just denote the time evolution of a measure \(\mu\) under the process dynamics:
\begin{tm}
In the sense of bounded test functions on \(\Omega\times\Omega\),
\be
\ba
\frac\di{\di t}(\tau_k\mush)S(t)&=\bigl(\e{\te(\vr+1)}-\e{\te(\vr)}\bigr)\cdot(\tau_{k+1}\mush-\mush)\\
&\quad+\bigl(\e{-\te(\vr)}-\e{-\te(\vr+1)}\bigr)\cdot(\tau_{k-1}\mush-\mush).
\ea\label{eq:srwgen}
\ee
\end{tm}
The first interesting consequence of this theorem is that the measure \(\mush\) on a coupled pair evolves into a linear combination of its shifted versions. Second, notice that \eqref{eq:srwgen} is the Kolmogorov equation for an asymmetric simple random walk. Indeed, this theorem implies the following
\begin{cor}\label{cr:rwtr}
Let the pair \((\un\xi^-(0),\,\un\xi(0))\) have initial distribution \(\mush\) defined by \eqref{eq:mushdf}. Then its later distribution evolves into a linear combination of translated versions of \(\mush\): at time \(t\) the pair \((\un\xi^-(t),\,\un\xi(0))\) has distribution
\[
\mush S(t)=\sum_{k=-\infty}^\infty P_k(t)\cdot\tau_k\mush,
\]
where \(P_k(t)\) is the transition probability at time \(t\) from the origin to \(k\) of a continuous time asymmetric simple random walk with jump rates
\[
\e{\te(\vr+1)}-\e{\te(\vr)}\text{ to the right and }\e{-\te(\vr)}-\e{-\te(\vr+1)}\text{ to the left.}
\]
In particular, \(Q^\xi(\cdot)\), started from an environment \(\mush\), is a continuous time asymmetric simple random walk with these rates.
\end{cor}
Although the corollary is quite natural, let us give a formal proof here. First some notation. \((\un\xi^-(\cdot),\,\un\xi(\cdot))\) will denote a pair of processes evolving under the basic coupling, \(g\) will be a bounded function on the path space of such a pair, and for shortness we introduce \(\Theta_t\) for the whole random path, shifted to time \(t\): \(\Theta_t=(\un\xi^-(t+\cdot),\,\un\xi(t+\cdot))\). Expectation of the process, started from \(\tau_k\mush\), will be denoted by \(\Ev^{(k)}\). Notice that under \(\Ev^{(k)}\) we a.s.\ have a single position \(Q^\xi(t)\) where the coupled pair differ by one, this is the position of the single conserved second class particle. With some abuse of notation we also use \(\Ev^{(\un\xi^-,\,\un\xi)}\) for the evolution of the pair \((\un\xi^-(\cdot),\,\un\xi(\cdot))\), started from the specific initial state \((\un\xi^-,\,\un\xi)\).

We aim for proving the semigroup property of \(S(\cdot)\). The first step is
\begin{lm}\label{lm:restart}
Given times \(0<s<t\) and \(k\in\Zb\),
\[
\Ev^{(0)}\bigl[g(\Theta_t)\,|\,Q^\xi(s)=k]=\Ev^{(k)}[g(\Theta_{t-s})].
\]
\end{lm}
\begin{proof}
The left hand-side is
\begin{multline*}
\frac{\Ev^{(0)}\bigl[g(\Theta_t)\,;\,Q^\xi(s)=k]}{\Pv^{(0)}\{Q^\xi(s)=k\}}\\
\ba
&=\frac{\Ev^{(0)}\bigl[\Ev^{(\un\xi^-(s),\,\un\xi(s))}g(\Theta_{t-s})\,;\,Q^\xi(s)=k]}{\Pv^{(0)}\{Q^\xi(s)=k\}}\\
&=\frac{\sum\limits_{j\in\Zb}\Pv^{(0)}\{Q^\xi(s)=j\}\Ev^{(j)}\bigl[\Ev^{(\un\xi^-(0),\,\un\xi(0))}g(\Theta_{t-s})\,;\,Q^\xi(0)=k]}{\Pv^{(0)}\{Q^\xi(s)=k\}}\\
&=\frac{\Pv^{(0)}\{Q^\xi(s)=k\}\Ev^{(k)}\bigl[\Ev^{(\un\xi^-(0),\,\un\xi(0))}g(\Theta_{t-s})\,;\,Q^\xi(0)=k]}{\Pv^{(0)}\{Q^\xi(s)=k\}}\\
&=\Ev^{(k)}[g(\Theta_{t-s})],
\ea
\end{multline*}
where in the second equality we used that the distribution at time \(s\) is a linear combination of shifted versions of \(\mush\).
\end{proof}
Next we prove the Markov property for \(Q^\xi(\cdot)\).
\begin{lm}
Let \(n>0\) be an integer, \(\vp_i\), \(i=0,\,\dots,\,n\) bounded functions on \(\Zb\), and \(0=t_0<t_1<\dots<t_n\). Then
\[
\Ev^{(0)}\prod_{i=1}^n\vp_i\bigl(Q^\xi(t_i)-Q^\xi(t_{i-1})\bigr)=\prod_{i=1}^n\Ev^{(0)}\vp_i\bigl(Q^\xi(t_i-t_{i-1})\bigr).
\]
\end{lm}
\begin{proof}
The statement is trivially true for \(n=1\). We proceed by induction, and assume the statement is true for \(n-1\). Then
\begin{multline*}
\Ev^{(0)}\prod_{i=1}^n\vp_i\bigl(Q^\xi(t_i)-Q^\xi(t_{i-1})\bigr)\\
\ba
&=\sum_{j\in\Zb}\Pv^{(0)}\{Q^\xi(t_1)=j\}\vp_1(j)\cdot\Ev^{(0)}\Bigl[\prod_{i=2}^n\vp_i\bigl(Q^\xi(t_i)-Q^\xi(t_{i-1})\bigr)\,|\,Q^\xi(t_1)=j\Bigr]\\
&=\sum_{j\in\Zb}\Pv^{(0)}\{Q^\xi(t_1)=j\}\vp_1(j)\cdot\Ev^{(j)}\prod_{i=2}^n\vp_i\bigl(Q^\xi(t_i-t_1)-Q^\xi(t_{i-1}-t_1)\bigr)\\
&=\sum_{j\in\Zb}\Pv^{(0)}\{Q^\xi(t_1)=j\}\vp_1(j)\cdot\Ev^{(0)}\prod_{i=2}^n\vp_i\bigl(Q^\xi(t_i-t_1)-Q^\xi(t_{i-1}-t_1)\bigr)\\
&=\sum_{j\in\Zb}\Pv^{(0)}\{Q^\xi(t_1)=j\}\vp_1(j)\cdot\prod_{i=2}^n\Ev^{(0)}\vp_i\bigl(Q^\xi(t_i-t_{i-1})\bigr)\\
&=\prod_{i=1}^n\Ev^{(0)}\vp_i\bigl(Q^\xi(t_i-t_{i-1})\bigr).
\ea
\end{multline*}
The second equality uses Lemma \ref{lm:restart}, the third one uses the fact that \(\phi\)'s only depend on \(Q^\xi\)-differences, and the fourth one follows from the induction hypothesis.
\end{proof}
\begin{proof}[Proof of Corollary \ref{cr:rwtr}]
We know that at any fixed time \(t>0\) the distribution of \((\un\xi^-(t),\,\un\xi(t))\) is a linear combination of shifted versions of \(\mush\). The shift is traced by the second class particle \(Q^\xi(t)\), therefore the differential equation
\begin{multline}
\frac\di{\di t}\Pv^{(0)}\{Q^\xi(t)=k\}\\
\ba
&=\bigl(\e{\te(\vr+1)}-\e{\te(\vr)}\bigr)\cdot\bigl(\Pv^{(0)}\{Q^\xi(t)=k+1\}-\Pv^{(0)}\{Q^\xi(t)=k\}\bigr)\\
&\quad+\bigl(\e{-\te(\vr)}-\e{-\te(\vr+1)}\bigr)\cdot\Bigl(\Pv^{(0)}\{Q^\xi(t)=k-1\}-\Pv^{(0)}\{Q^\xi(t)=k\}\Bigr)
\ea\label{eq:qsrwgen}
\end{multline}
follows from \eqref{eq:srwgen}. In the above lemmas, we also proved that \(Q^\xi(t)\) is Markovian (annealed w.r.t.\ the initial distribution of \((\un\xi^-,\,\un\xi)\)). As there exists only one Markovian process with Kolmogorov equation \eqref{eq:qsrwgen} of the simple asymmetric random walk, we conclude that the process \(Q^\xi(\cdot)\) with initial environment \(\mush\) is an asymmetric simple random walk with rates as stated in the Corollary.
\end{proof}
\begin{lm}
Let $(\uomm,\,\uomp)$ be a pair of processes in basic coupling, started from distribution \eqref{eq:umhdef}, with second class particle $Q(t)$. Then there exist constants \(0<\al_0,\,C<\infty\) such that 
\[
\Pv\{|Q(t)|>K\}\le\e{-CK}
\]
whenever \(K>\al_0 t\) and $t$ is large enough.
\end{lm}
Notice that this implies that Assumption \ref{ass:main} holds for the TAEBLP.
\begin{proof}
The proof uses auxiliary processes to connect the above arguments to the setting of Assumption \ref{ass:main}. Define the pair
\[
(\la_i,\,\vr_i):\,=\left\{
\ba
&(\vr,\,\vr+1),&&\text{for }i\le0,\\
&(\vr,\,\vr),&&\text{for }i>0.
\ea
\right.
\]
Draw the pair \((\un\ze(0),\,\un\xi(0))\) from the product distribution of coupling measures \eqref{eq:cmu}
\[
\bigotimes_{i\in\Zb}\mu^{\la_i,\vr_i}.
\]
Then \(\un\xi(0)\) has distribution
\[
\bigotimes_{i\le0}\mu^{\vr+1}\cdot\bigotimes_{i<0}\mu^\vr,
\]
in agreement with \eqref{eq:mushsec}.

Let now the pair \((\un\ze(\cdot),\,\un\xi(\cdot))\) evolve in the basic coupling, and let them play the role of \((\un\eta(\cdot),\,\un\om(\cdot))\) of Section \ref{sc:mcbl}. This results in the pair \((\un\ze(\cdot),\,\un\ze^+(\cdot))\) with a second class particle \(Q^\ze(\cdot)\) and the pair \((\un\xi^-(\cdot),\,\un\xi(\cdot))\) with a second class particle \(Q^\xi(\cdot)\) such that \(Q^\ze(t)\le Q^\xi(t)\), see Lemma \ref{lm:zryz}. Therefore the random walk result in Corollary \ref{cr:rwtr} on \(Q^\xi(\cdot)\) yields the desired estimate for \(Q^\ze(t)\). 
Finally, notice that the distribution of \(\un\om^-(0)\) in Assumption \ref{ass:main} and of \(\un\ze(0)\) above only differ by \(\om^-_0(0)\sim\wih\mu^\vr\), while \(\ze_0(0)\sim\mu^\vr\). Therefore
\[
\ba
\Pv\{Q(t)>K\}&=\sum_{z=-\infty}^\infty\Pv\{Q(t)>K\,|\,\om^-_0(0)=z\}\cdot\mu^\vr(z)^\frac12\Bigl(\frac{\wih\mu^\vr(z)^2}{\mu^\vr(z)}\Bigr)^\frac12\\
&=\sum_{z=-\infty}^\infty\Pv\{Q^\ze(t)>K\,|\,\ze_0(0)=z\}\cdot\mu^\vr(z)^\frac12\Bigl(\frac{\wih\mu^\vr(z)^2}{\mu^\vr(z)}\Bigr)^\frac12\\
&\le\Bigl[\sum_{z=-\infty}^\infty\Pv\{Q^\ze(t)>K\,|\,\ze_0(0)=z\}\cdot\mu^\vr(z)\Bigr]^\frac12\cdot\Bigl[\sum_{y=-\infty}^\infty\frac{\wih\mu^\vr(y)^2}{\mu^\vr(y)}\Bigr]^\frac12\\
&=\Pv\{Q^\ze(t)>K\}^\frac12\cdot\Bigl[\sum_{y=-\infty}^\infty\frac{\wih\mu^\vr(y)^2}{\mu^\vr(y)}\Bigr]^\frac12.
\ea
\]
We are done as soon as we show that \(\wih\mu^\vr(y)/\mu^\vr(y)\) is uniformly bounded in \(y\). With the exponential rates \eqref{eq:fdef} one obtains from \eqref{eq:muhat}
\[
\frac{\wih\mu^\vr(y)}{\mu^\vr(y)}=C\sum_{z=y+1}^\infty(z-\vr)\e{-\frac\beta2(z-\frac\te\beta)^2+\frac\beta2(y-\frac\te\beta)^2}=C\sum_{k=1}^\infty(k+y-\vr)\e{-\frac\beta2k^2-\beta ky+\te k}.
\]
This is uniformly bounded for large \(y\)'s since then \(y\e{-\beta y}<1\). For large negative \(y\)'s one uses the equivalent form
\[
\wih\mu^\vr(y):\,=\frac{1}{\Vv^\vr(\om_0)}\sum_{z=-\infty}^y(\vr-z)\mu^\vr(z)
\]
of \eqref{eq:muhat} and writes
\[
\frac{\wih\mu^\vr(y)}{\mu^\vr(y)}=C\sum_{z=-\infty}^y(\vr-z)\e{-\frac\beta2(z-\frac\te\beta)^2+\frac\beta2(y-\frac\te\beta)^2}=C\sum_{k=1}^\infty(k-y+\vr)\e{-\frac\beta2k^2+\beta ky-\te k}
\]
which is again uniformly bounded for large negative \(y\) values.

To show a lower bound on $Q(t)$, start with
\[
 (\la_i,\,\vr_i):\,=\left\{
  \ba
   &(\vr,\,\vr),&&\text{for }i<0,\\
   &(\vr-1,\,\vr),&&\text{for }i=0,\\
   &(\vr,\,\vr-1),&&\text{for }i>0,
  \ea
 \right.
\]
and the coupled pair \((\un\ze(0),\,\un\xi(0))\) in distribution
\[
\bigotimes_{i\in\Zb}\mu^{\la_i,\vr_i}.
\]
Now the roles of the pair \((\un\ze(\cdot),\,\un\ze^+(\cdot))\) with a second class particle \(Q^\ze(\cdot)\) and the pair \((\un\xi^-(\cdot),\,\un\xi(\cdot))\) with a second class particle \(Q^\xi(\cdot)\) are interchanged and we have \(Q^\ze(t)\ge Q^\xi(t)\). The random walk estimate on $Q^\xi$ and a Radon-Nikodym estimate similar to the one above completes the proof of the lower bound.
\end{proof}

\appendix
\appendixpage

\section{Covariance identities for bricklayer process with exponential rates}
The purpose of this appendix is to prove
the variance formula for
stationary BLP  under the following exponential bound
assumption on rates: for some $C,\beta<\infty$, 
\[
\rate(\om_0)\le Ce^{\beta\abs{\om_0}}.
\]
 Assume the height process
is normalized initially by $h_0(0)=0$. 

\begin{tm} 
Fix $z\in\Zb$. 
 In the stationary  infinite volume process with time-marginal
distribution $\un\om(t)\sim\un\mu^\te$, 
\be
\Vv[h_z(t)]=\sum_{n\in\Zb} \abs{n-z} \Cov[\om_n(t),\om_0(0)].
\label{varh1}\ee
\label{varh1thm}\end{tm}

Formula \eqref{varh1}  was already proved
in  \cite{varj2nd} for a general class of processes
that includes ZRP and BLP.  However, 
 the proofs in \cite{varj2nd} were carried 
out under the assumption that certain semigroup and generator 
calculations work.  This presents no problem
when the single-site state space
is compact (such as
exclusion processes
) for then one has a strongly
continuous semigroup on the Banach space of
continuous functions on the state space of the
process.  For BLP with superlinear rates, only
some rudimentary features of the usual semigroup
picture have been established in \cite{exists}.
Hence the need to  justify \eqref{varh1}.  
We use the finite-volume auxiliary processes as introduced in \cite{exists}.

To prove Theorem \ref{varh1thm} we show that the
infinite-volume stationary process is a limit of
finite-volume $\lrt$ processes as $-\ell,\err\to\infty$.
A preliminary form of \eqref{varh1} is 
true for an $\lrt$ process by simple
counting. (See \eqref{varhaux7} below and its 
expanded form on lines \eqref{varhaux10a}--\eqref{varhaux10d}.)
  The technical work goes into establishing
moment bounds that are uniform over $\ell<0<\err$.
These in turn permit us to take the 
$-\ell,\err\to\infty$  limit 
in the proto-formula  \eqref{varhaux10a}--\eqref{varhaux10d}.  

This appendix is based on the construction of the 
infinite-volume BLP $\hhgt(t)$ given in \cite{exists}.
Article \cite{exists} utilized two types
of finite-volume processes:  the  {\sl $\lr$ processes}
whose height variables were denoted by  $\hhgt^{\lr}(t)$, 
and the {\sl $\lrt$ processes} with height variables 
denoted by  $\ghgt^{\lrt}(t)$.  

The $\lr$ evolution
is a straightforward restriction of the full system into 
the finite interval $\lr$, with generator 
\[
L^{\lr}\vp(\un\om)=\sum_{i=\ell}^{\err-1}\left[f(\om_i)+f(-\om_{i+1})\right]\cdot\left[\vp(\un\om^{(i,\,i+1)})-\vp(\un\om)\right].
\]
This generator defines a countable state space Markov process 
that evolves over the sites  $\ell,\,\dots,\,\err$: the jump $\un\om\to\un\om^{(i,\,i+1)}$ happens at rate $f(\om_i)+f(-\om_{i+1})$, independently at different sites $i$, but only for $\ell\leq i\leq \err-1$. Columns outside
the interval $[\ell, \err]$ are frozen for all time.  
The virtue of this process is monotone dependence on the 
interval $\lr$.    In \cite{exists} the infinite-volume process $\hhgt(\cdot)$ was  
defined as the a.s.\ increasing limit
of the processes  $\hhgt^{\lr}(t)$ as $-\ell,\err\to\infty$.

  The  $\lrt$ process has also
the correct boundary currents that make the i.i.d.\ product
measures $\un\mu^\te$ invariant for the finite collection
of increment variables $\{\om^{\lrt}_i: \ell\le i\le\err\}$.  The generator is 
\[
\ba
G^{\lrt}\vp(\un\om)&=\sum_{i=\ell}^{\err-1}\,[f(\om_i)+f(-\om_{i+1})]\cdot\left[\vp(\un\om^{(i,\,i+1)})-\vp(\un\om)\right]\\
+\,&[\e{\te}+f(-\om_{\ell})]\cdot[\vp(\un\om^{(\ell-1,\,\ell)})-\vp(\un\om)]\\
\,+\,&[\e{-\te}+f(\om_\err)]\cdot[\vp(\un\om^{(\err,\,\err+1)})-\vp(\un\om)].
\ea
\]

For a concrete construction of the processes,
we imagine that bricklayer at site $i$ has two
unit rate  Poisson
processes $N_i^{(L)}$ and $N_i^{(R)}$
on the first quadrant $\Rb_+^2$ of the plane. 
These govern his brick-laying
action to the left $(L)$ and right $(R)$. 
A  Poisson point $(t,y)$ in $N_i^{(L)}$ 
such that $y\le \rate(-\om_i(t-))$ signals a brick to be laid
on $[i-1,i]$ at time $t$, 
while a  Poisson point $(t,y)$ in $N_i^{(R)}$ 
such that $y\le \rate(\om_i(t-))$ signals a brick to be laid
on $[i,i+1]$ at time $t$.
Shift of these planar Poisson processes by time \(t\) will be denoted by \(S_tN\).

The construction of the finite
systems in terms of these Poisson processes
provides the usual jump chain-holding time
construction of a continuous time Markov chain on a
countable state space.  After each jump, the holding
time and the next state are read from the Poisson
processes.  By the strong Markov property of the
Poisson processes this is equivalent to looking
at freshly sampled exponential variables with  
appropriate rates. In both   $\lr$ processes 
and  $\lrt$ processes
the increment
$H(t)-H(0)$  of the maximal height 
\be
H(t)=\max_{\ell-1\leq j\leq \err}h_j(t)
\label{defHmax}\ee
 is bounded
by a Poisson process.  Hence  
explosions do not happen \cite[Sect.~3.1-3.2]{exists}. 

In the next section we show that the $\lrt$ processes
converge to the infinite volume stationary process.
Then we develop moment bounds through martingales,
uniformly in $\ell<0<\err$. 
 After Section
\ref{sc:lrtconv} we drop the superscripts $\lr$ and $\lrt$
to ease notation.

\subsection{Convergence of $\lrt$ processes}
\label{sc:lrtconv}
The infinite-volume process $\hhgt(\cdot)$ is 
defined as the a.s.\ increasing limit
of the processes  $\hhgt^{\lr}(\cdot)$
\cite[Theorem 2.2]{exists}.   
 Lemma 7.1 in \cite{exists}
shows that if the initial state $\un\om(0)$ is 
$\un\mu^\theta$-distributed
then the increment process $\un\om(\cdot)$  is stationary. 

  Calculations in this appendix are done in
a  stationary finite-volume $\lrt$ process.
Hence we need to show that this process also converges as
$-\ell,\err\to\infty$ to the stationary infinite-volume
process.  

It will be 
 convenient to represent the processes by measurable
mappings of  the initial
configuration $\hhgt$ and the Poisson clocks $\Nconf$ on $(0,t]$:
\[\ghgt^{\lrt}(t)=\Psi_t^{\lrt}(\hhgt,\Nconf),\ 
\hhgt^{\lr}(t)=\Phi_t^{\lr}(\hhgt,\Nconf)  
\quad\text{and}\quad
\hhgt(t)=\Phi_t(\hhgt,\Nconf). \]
Then the construction of the process $\hhgt(t)$ in \cite{exists} 
can be expressed as follows: 
for any initial $\hhgt\in\wt\Omega$, any $m,T<\infty$, 
\be
\Phi_{t,i}(\hhgt,\Nconf)=
 \Phi_{t,i}^{\lr}(\hhgt,\Nconf)
\quad\text{for  $-m\le i\le m$ and $0\le t\le T$}
\label{Ahlim} \ee
for all large enough $-\ell,\err$, almost surely.

\begin{lm}  Let $\hhgt(t)$ be the infinite-volume
process whose increment  process $\un\om(t)$ is stationary
with marginal $\un\om(t)\sim\un\mu^\theta$. 
As $-\ell,\err\to\infty$, $\ghgt^{\lrt}(\cdot)\to 
 \hhgt(\cdot)$ in the following sense:  given 
 any $m,T<\infty$,   
\[\text{$g^{\lrt}_i(t)=h_i(t)$
for  $-m\le i\le m$ and $0\le t\le T$}\]
for all large enough $-\ell,\err$.
\label{lrtconvlm}\end{lm}

\begin{proof}  
The proof of Lemma 7.1 on p.~1243 of \cite{exists} shows that
there exists a (nonrandom) time $t_0=t_0(\te)>0$ such that, for any $m$, 
$g^{\lrt}_i(t)=h^{\lr}_i(t)$ for $-m\le i\le m$ and $0\le t\le t_0$
if $-\ell$ and $\err$ are large enough. Combined
with \eqref{Ahlim}  we have the statement on the time interval 
$[0,t_0]$. 

To get to an arbitrary time $T$, 
 the claim is proved by induction up to time $kt_0>T$, for a
positive integer $k$.   Since   process $\hhgt(\cdot)$ stays
in $\wt\Omega$,  \eqref{Ahlim}  can be applied to restarted
processes.  

 The other ingredient of the induction step is
control of discrepancies between processes 
$g^{\lrt}_i(\cdot)$ and $h^{\lr}_i(\cdot)$. 
For this purpose let $\Bc$ be the event that 
in the process $\{\ghgt^{\abt}(kt_0+t): 0\le t\le t_0\}$ either 
all columns
in the range $\{\fl{3\ell/4}, \dotsc, \fl{\ell/2}+1\}$  grew, or 
 all  columns in the range $\{\fl{\err/2},\dotsc,\ce{3\err/4}\}$ grew. 
By \cite[Corollary 5.5]{exists},
 $\Pv(\Bc)\le C_1e^{-C_2(\abs{\ell}\wedge\err)} $, 
and this bound is independent of $a,b$.  

Discrepancies
originate at the  edges $\ell$ and $\err$.  
On the event  $\Bc^c$ the   $\abt$ and $\lr$ processes started from a common initial
configuration are indistinguishable inside $(\ell/2,\err/2)$ because 
 discrepancies have not had a chance to penetrate into  $(\ell/2,\err/2)$. 
By the monotonicity properties of these processes
the columns of the $\lrt$ process never go below those of 
the $\lr$ process \cite[Lemma 3.2]{exists}.  Hence it is enough
to use the event $\Bc^c$ 
to suppress growth in the $\lrt$ process.   We leave the details
of this induction argument to the reader. 
\end{proof}

\subsection{Martingales}
Since the rates are unbounded, we begin by stating  a 
  general lemma about countable Markov chains.
Let $S$ be a countable state space, $Q$ a
generator matrix,  $q_x=-q_{x,x}=\sum_{y:y\ne x} q_{x,y}$
the total rate to jump away from state $x\in S$.
Let $P(t)$ be the semigroup
associated to $Q$, in other words the minimal positive
solution of the backward equation $P'(t)=QP$, $P(0)=I$
\cite{norris}. 
The next lemma is proved with standard techniques.  
 
\begin{lm}  Let $\nu$ be an initial distribution on $S$ and $T<\infty$.  
Assume 
\be
\int_0^T E^\nu[q_{X_s}]\di s
= \int_0^T \sum_x\nu(x)\sum_y p_{x,y}(s)q_y\di s<\infty.
\label{fwass}\ee
Let  $\varphi$ be a  bounded function on $S$.
Then  under $P^\nu$,  for $t\in[0,T]$, 
the process 
$\sigma(t)=\int_0^t Q\varphi(X_s)\di s $
is well-defined and in $L^1(P^\nu)$, 
and the following process is a martingale: 
\be
M_t= \varphi(X_t)-\varphi(X_0)-\int_0^t Q\varphi(X_s)\di s . 
\label{mg1}\ee
\end{lm}

Now apply this to the $\lr$ and  $\lrt$ processes. To simplify notation, let
\[
f_i(\un\om(s))=f(\om_i(s))+f(-\om_{i+1}(s))
\]
be the rate of growth of the column height \(h_i(s)\).

\begin{lm} Fix $\ell<0<\err$. Consider 
either the $\lr$ process or the  $\lrt$ process
and in either case denote the height variables 
by $h_i(t)$. 
Let $\nu$ be an initial distribution 
 such that for some $c>\beta$ 
\be E^\nu \bigl( e^{c\abs{h_0}\,+\,
c\sum_{i=\ell}^{\err}\abs{\om_i}}\bigr)
<\infty. \label{hmgass1}\ee
Then for any $1\le p<\infty$ 
and any index  $i$  this process is a  martingale: 
\be
M_t=h^p_i(t)-h^p_i(0) -\int_0^t 
\rate_i(\un\om(s))\bigl( (h_i(s)+1)^p-h^p_i(s)\bigr)\di s.
\label{hlrmg}\ee
\end{lm}

\begin{proof}
In the $\lr$ process the total jump rate
out of state $\hhgt$ is 
\[
q_{\hhgt} =\sum_{i=\ell}^{\err-1} \rate_i(\un\om) 
\le 2\sum_{i=\ell}^{\err} e^{\beta\abs{\om_i}}
\le 2(\err-\ell+1) e^{\beta\sum_{i=\ell}^{\err}\abs{\om_i}}.
\]
In the $\lrt$ process  the total jump rate
out of state $\ghgt$  is
\[
q^\te_{\ghgt}=q_{\ghgt}+ \rate(-\om_\ell)+\rate(\om_\err)
+e^\te+e^{-\te}
\le 2(\err-\ell+1) e^{\beta\sum_{i=\ell}^{\err}\abs{\om_i}}
+e^\te+e^{-\te}.
\]

In the $\lr$ process
under a fixed initial configuration $\hhgt$,  
\[
\Ev^{\hhgt}\Bigl[\,\sup_{s\in[0,T]}
e^{\beta\sum_{i=\ell}^{\err}\abs{\om_i(s)}} \,\Bigr]
\le
e^{\beta\sum_{i=\ell}^{\err}\abs{\om_i}} \Ev[\,
e^{\beta Y(T)}\,]
\]
where $Y(\cdot)$ is a Poisson process of rate 
$\lambda=2\rate(0)(\err-\ell)$. This comes from the observation that 
the process
\[v(t)= \sum_{i=\ell}^{\err}\lvert \om_i(t)\rvert
\,-\, \sum_{i=\ell}^{\err}\lvert \om_i(0)\rvert
\]
increases only when a local maximum column grows, and
this happens at rate at most $2\rate(0)(\err-\ell)$. 
Then, under the  initial distribution $\nu$,
\be\begin{split}
\Ev^\nu\Bigl[\,\sup_{s\in[0,T]} q_{\hhgt(s)}\,\Bigr] 
&\le \Ev^\nu\Bigl[\,\sup_{s\in[0,T]}
e^{\beta\sum_{i=\ell}^{\err}\abs{\om_i(s)}} \,\Bigr]\\
&\le
\Ev^\nu\bigl( e^{\beta\sum_{i=\ell}^{\err}\abs{\om_i}}\bigr)  
\cdot\exp\{ \lambda(e^{\beta T}-1)\} 
<\infty. \end{split}
\label{auxexpbound1}\ee

In the $\lrt$ process the same idea works
except the rate for the process that bounds $v(t)$ is 
$\lambda^\te=2\rate(0)(\err-\ell+1)  +e^\te+e^{-\te}$.

Let $F(x)=(b\wedge x)\vee(-b)$ be a truncation function.
Now that assumption \eqref{fwass} has been verified, 
for any integer $0<b<\infty$ and $\ell\le i\le \err-1$
 \eqref{mg1} gives
 the martingale
\[
M^{(b)}_t=F(h_i(t))^p\;-\;F(h_i(0))^p\; -\;\int_0^t 
\rate_i(\un\om(s)) 
\bigl( F(h_i(s)+1)^p-F(h_i(s))^p\bigr) \di s.
\]
So for an event  $A\in\Fc_s$ and  $s<t$ 
\begin{align*}
&\Ev^\nu\bigl[F(h_i(t))^p \ind_A\bigr]
-\Ev^\nu\bigl[F(h_i(s))^p  \ind_A\bigr]\\
 &\qquad = \int_s^t 
\Ev^\nu\bigl[\rate_i(\un\om(u)) 
\bigl( F(h_i(u)+1)^p-F(h_i(u))^p\bigr) \ind_A\bigr]
 \di u.
\end{align*}
As $b\nearrow\infty$ dominated convergence takes each term
to the desired limit. This is justified by the following.
Restrict to an interval $s,t\in [0,T]$. 
The rate in the last expectation is bounded as in
\[\rate_i(\un\om(u))\le \sup_{s\in[0,T]} 
e^{\beta(\abs{\om_{i}(s)}+\abs{\om_{i+1}(s)}) }
\]
and the random variable on the right has a finite $L^p$-norm
for some $p>1$ by a  bound of the type  in \eqref{auxexpbound1}
and the assumption that $c>\beta$ in \eqref{hmgass1}.
For the height we have
\[\abs{h_i(t)}=h_i^+(t)+h_i^-(t)\le H(t)+h_i^-(0) \le 
(H(t)-H(0))  +H(0)+ \abs{h_i(0)}
\]
where $H(t)$ is the maximal height of \eqref{defHmax}.
The increment $H(t)-H(0)$ is
stochastically bounded by a Poisson process  
while $H(0)$ and $\abs{h_i(0)}$
 have all moments by assumption \eqref{hmgass1}.
 We conclude that  \eqref{hlrmg} is a 
martingale.
\end{proof} 

\subsection{Bounds for the $\lrt$ process}
Henceforth consider stationary  $\lrt$ processes with
$\mu^\te$ marginals and 
 initial height 
normalized by $h_0(0)=0$.   The increment process
is denoted by $\un\om(t)$ and the height process by
$\hhgt(t)$. 
In this process the column heights $h_{i}(t)$ 
for $i\le\ell-2$ and $i\ge\err+1$ 
are frozen at their initial values.  We start
with moment bounds that hold uniformly in $\ell<0<\err$.

\begin{lm}  Fix $\te$.  
For $1\le p<\infty$ there is a constant $C=C(p,\te)$
such that 
 in all  stationary 
$\lrt$ processes with marginal distribution $\mu^\te$,
 for all $t\ge 0$ and $i\in\Zb$,  
\[
\Ev\bigl[(h_i(t)-h_i(0))^p\,\bigr] \le e^{Ct}. 
\] 
The bound is valid also for the boundary columns 
$i=\ell-1$ and $i= \err$,  and also for the
infinite-volume stationary process with marginal distribution
$\mu^\te$. 
\label{hmombdlm}\end{lm}

\begin{proof}  Abbreviate
$\bar h_i(t)=h_i(t)-h_i(0)$.
It suffices to consider the $\lrt$ processes because the bound
extends to the $-\ell,\err\to\infty$ limit by Fatou's lemma
and Lemma \ref{lrtconvlm}. 
The increments $\om_i$ have all exponential moments under 
$\mu^\te$ so assumption \eqref{hmgass1} is satisfied. 
In particular, $\bar h_i(t)$ has all moments. 

Fix $\ell<0<\err$.  Columns $h_i$ for $i\notin[\ell-1,\err]$ 
are frozen at their initial values and need no argument. 
Consider the case $\ell\le i\le \err-1$. 
It suffices to consider a positive integer $p$. 
  In the next calculation,
use martingale  \eqref{hlrmg}, the fact that 
$\bar h_i(s)$ is a nonnegative integer to bound a lower power
with a higher one, and apply H\"older's inequality:
\begin{align*}
\Ev\bar h_i(t)^p &= \Ev\int_0^t  \rate_i(\un\om(s))\sum_{k=0}^{p-1}\binom{p}{k}\bar h_i^k(s)\di s\\
&\le C_p\int_0^t  \Ev\bigl[ \rate_i(\un\om(s)) \bigl(1\vee\bar h_i(s)\bigr)^{p-1}\bigr]\di s \\
&\le C_{p,\te} \int_0^t \bigl(\Ev[\rate_i(\un\om(s))^{p}]\bigr)^{1/p} \bigl(\Ev[1\vee\bar h_i(s))^p]\bigr)^{(p-1)/p}\di s\\
&\le C_{p,\te} \int_0^t \bigl(\Ev\bar h_i(s)^p+1\bigr) \di s.
\end{align*}

Rewrite this as 
\[
\Ev\bar h_i(t)^p+1 \le 1+ 
C_{p,\te} \int_0^t \bigl(\Ev\bar h_i^p(s)+1\bigr) \di s
\]
and now Gronwall's inequality gives the conclusion.

The boundary columns  $i=\ell-1$ and $i=\err$ are handled
similarly, with the only difference that the rates include
also the constant terms $e^\te$ or $e^{-\te}$.  
\end{proof}

Fix a path $z(t)$ in $\Zb$ with $z(0)=0$ and
since  $t$ is fixed abbreviate  $z=z(t)$.   
 Define 
\[
\pcount_+(t)=\sum_{n= z(t)+1}^{\err+1}\om_n(t)
\quad\text{and}\quad 
\pcount_-(t)=\sum_{n=\ell-1}^{z(t)}\om_n(t).
\]
Due to the frozen columns at the boundaries and the normalization $h_0(0)=0$
 we have the identity
\[
h_z(t)=\pcount_+(t)-\pcount_+(0)=-\pcount_-(t)+\pcount_-(0)
\]
from which 
\be
\Vv[h_z(t)]
= \Cov\bigl[ -\pcount_-(t)+\pcount_-(0)\,,\,
\pcount_+(t)-\pcount_+(0)\bigr].
\label{varhaux7}\ee
For each time 
$t$  the increment variables $\{\om_i(t): \ell\le i\le\err\}$
are i.i.d.\ $\mu^\te$-distributed. 
The process is independent of the initial
values $\om_{\ell-1}(0)$ and $\om_{\err+1}(0)$ of the 
boundary increments. This is   because the 
bricklayer at site $\ell-1$ lays bricks to his right at
rate $e^\te$ regardless of the value of the increment at site $\ell-1$,
and similarly the left action of the bricklayer
 at site $\err+1$ has constant rate $e^{-\te}$.  But the time increment
\be
\om_{\ell-1}(t)-\om_{\ell-1}(0)
= -h_{\ell-1}(t)+h_{\ell-1}(0)
\label{omellincr}\ee
and the same for $\om_{\err+1}(t)$ are needed  for $h_z(t)$.
Expanding and removing vanishing covariances from \eqref{varhaux7}
leaves
\begin{align}  
&\Vv[h_z(t)] 
=\sum_{\substack{\ell\le i\le0\\z<j\le\err}} 
\Cov(\om_i(0),\om_j(t)) 
\;+\; 
\sum_{\substack{\ell\le i\le z\\0<j\le\err}} 
\Cov(\om_i(t),\om_j(0))\label{varhaux10a}\\
&\qquad \;+\;
\sum_{\ell\le i\le0} 
\Cov(\om_i(0),\om_{\err+1}(t)) 
\;+\; 
\sum_{0<j\le\err} 
\Cov(\om_{\ell-1}(t),\om_j(0))\label{varhaux10b}\\
&\qquad \;-\;
\sum_{\ell\le i\le z} 
\Cov(\om_i(t),\om_{\err+1}(t)) 
\;-\; 
\sum_{z<j\le\err} 
\Cov(\om_{\ell-1}(t),\om_j(t))\label{varhaux10c}\\
&\qquad \;-\;  \Cov(\om_{\ell-1}(t),\om_{\err+1}(t)).
\label{varhaux10d}
\end{align}

We argue that the sums on lines 
\eqref{varhaux10b}--\eqref{varhaux10d} vanish
as $-\ell,\err\to\infty$ by showing that the covariances
decay exponentially in the spatial distance.

We illustrate with a term
$\Cov(\om_{i_0}(t),\om_{\err+1}(t))$
for a fixed $i_0\le z$
 from the first sum on line \eqref{varhaux10c}.
Consider $-\ell,\err$ large so that $\ell<i_0<z<\err$.  
 Let $w=\fl{(i_0+\err)/2}$ be a lattice point at or next
to the midpoint
between $i_0$ and $\err$.   
Let $\un\ze(t)$ and $\un\xi(t)$ be two further processes 
whose initial configurations agree with those of
 $\un\om(t)$, both for the increments
 and for the heights:  
\be \un\ze(0)=\un\xi(0)=\un\om(0) 
\quad\text{and}\quad
\hhgt^\ze(0)=\hhgt^\xi(0)=\hhgt(0). 
\label{initzetaxi}\ee

Processes  $\un\ze(t)$ and $\un\xi(t)$  follow the 
$\lrt$ evolution, except that 
the columns $\{\hhgt^\ze_i(t):i\ge w-1\}$ 
 with bases  in  $[w-1,\infty)$ are 
not permitted to grow, 
and the columns   $\{\hhgt^\xi_i(t):i\le w\}$ based in 
 $(-\infty,w]$ are similarly frozen. (Equivalently, replace the Poisson
clocks of the corresponding bricklayers with empty point measures.)  
To determine their dynamics, in addition to disjoint collections
of Poisson clocks,  $\un\ze(t)$ requires initial
increments $\{\om_i(0):i\le w-1\}$ while $\un\xi(t)$  requires initial
increments $\{\om_i(0):i\ge w+1\}$.
Thus $\{\ze_i(t):i\le w-1\}$ and 
$\{\xi_i(t):i\ge w+1\}$ are independent processes.
The height processes 
$\{h^\ze_i(t):i\le w-2\}$ and 
$\{h^\xi_i(t):i\ge w+1\}$ are not independent.
For example if $z>0$ 
 the initial heights $\{h^\xi_i(0):i\ge w+1\}$ are
very much influenced  by the initial increments $\{\ze_i(0):0<i\le z\}$.

On the sites where  $\un\ze(t)$ and $\un\xi(t)$ lay bricks 
they obey the same realizations of Poisson clocks as the
$\lrt$ process $\un\om(t)$. 
This coupling preserves the inequalities
\be
\text{$\hhgt^\ze(t)\le\hhgt(t)$ and $\hhgt^\xi(t)\le\hhgt(t)$.}
\label{coupling3}\ee

Let
\be
A=\{ \text{$\om_i(t)=\ze_i(t)$ for $i\le i_0$}\}
\quad\text{and}\quad  
B=\{\text{$\om_{\err+1}(t)=\xi_{\err+1}(t)$}\}.
\label{defAB4}\ee

Let us say that the space-time window $[a,b]\times(s_0,s_1]$
 {\sl is a block}  if $h_i(s_1)=h_i(s_0)$ for some 
$\ce{a}\le i\le \fl{b}-1$,
in other words, if some column failed to grow inside
 the space interval $[a,b]$ 
during time interval  $(s_0,s_1]$. The space-time window is not
a block if every column inside grew during $(s_0,s_1]$. 
The sense of the terminology is that a block does not permit
a discrepancy to pass.  
 For convenience
we do not require $a,b$ to be integers.   We restate 
Corollary 5.5 from \cite{exists}.  This lemma is valid
in the stationary process because there is control over the
spatial averages of the increments $\om_i$ and thereby
 control  over rates. 

\begin{lm} {\rm\cite[Cor.~5.5]{exists}}
There exist constants $k_0<\infty$ and $t_0>0$
and a function $0<G(s)<\infty$ for $s\in(0,t_0)$ such that 
$G(s)\nearrow\infty$ as $s\searrow 0$, and this bound holds:
if $s_1-s_0\le t_0$ and $b-a\ge k_0$ then 
\[
\Pv\{\text{$[a,b]\times(s_0,s_1]$
is not a block} \} \le e^{-(b-a)\cdot G(s_1-s_0)}. 
\]
For a given $\te$ this bound is valid for all $\lrt$ processes.
\end{lm}

Recalling that $t$ is fixed in the present discussion,
fix a positive integer $m$ and real $\tau>0$   so that
\[
\text{$t=m\tau$ and $\tau\in(0, t_0]$.}\]
   
On the event $A^c$ there must exist a sequence of 
times $0<t_{w-2}<t_{w-3}<t_{w-4}<\dotsm <t_{i_0}\le t$ such that
column $h_i$  grew at time $t_i$, $i_0\le i\le w-2$. 
This results from the 
 observation that in the range $i\le w-2$ 
 the first  discrepancy $h_i-h^\ze_i>0$
 in column heights at $i$ can appear only next to an existing
discrepancy.  Thus the leftmost discrepancy
$Q(t)=\min\{i:  h_i(t)-h^\ze_i(t)>0\}$ starts from the 
value $Q(0)=\infty$ due to \eqref{initzetaxi}, 
arrives at the boundary $w-1$ of the frozen region of 
$\hhgt^\ze$ at time $t_{w-1}$ when $h_{w-1}$ first grows,  
 and then moves to the left
one step at a time and always with a jump of the column
$h_{Q(t)-1}(t)$ that is not matched by a jump in
$h_{Q(t)-1}^\ze(t)$.  

Consequently event $A^c$ implies 
that at least one of the windows  
\[ \Bigl[w-2-j\frac{w-2-i_0}m, w-2-(j-1)\frac{w-2-i_0}m\Bigr]
\times[(j-1)\tau, j\tau], \quad 1\le j\le m,\]
 is
not a block.  For if all these windows were blocks, the sequence
of growths over intervals  $[w-i,w-i+1]$ at times $t_{w-i}$,
 $2\le i\le w-i_0$,
could not happen. This gives the bound
\[
\Pv(A^c)\le me^{-C(w-2-i_0)/m}\le me^{-C(\err-i_0)/m}
\]
where for the second inequality we 
  assumed that $\err$ is large enough relative to $z$.

The same argument for the propagation of discrepancies can be
repeated for event $B$ in \eqref{defAB4} to improve the bound to
\be
\Pv(A^c\cup B^c)\le 2me^{-C(\err-i_0)/m} \le C_1e^{-C_2(\err-i_0)}
\label{PABprob1}\ee  
 where we subsumed $m$ in the constants.  

Next we turn these
bounds into bounds on covariances. 
\be\begin{split}
&\Cov(\om_{i_0}(t),\om_{\err+1}(t)) 
=\Cov(\ze_{i_0}(t),\xi_{\err+1}(t)) \\
&\qquad\qquad +\Ev\bigl[(\om_{i_0}(t)\om_{\err+1}(t)
-\ze_{i_0}(t)\xi_{\err+1}(t))
\ind_{A^c\cup B^c}\bigr]\\
&\qquad\qquad -\Ev\bigl[(\om_{i_0}(t)-\ze_{i_0}(t))\ind_{A^c}\bigr]
\cdot\Ev\om_{\err+1}(t) \\
&\qquad\qquad\;-\; \Ev\bigl[(\om_{\err+1}(t)-\xi_{\err+1}(t))
\ind_{B^c}\bigr]\cdot\Ev\ze_{i_0}(t).
\end{split}\label{covbd3}\ee
By independence of  $\ze_{i_0}(t)$ and $\xi_{\err+1}(t)$
 the first covariance after the equality
sign vanishes.  
We claim that for $\ell\le i\le \err$,
\be
\Ev\abs{\om_{i}(t)}^p \;+\; 
\Ev\abs{\ze_{i}(t)}^p \;+\;
\Ev\abs{\om_{\err+1}(t)}^p  \;+\;
\Ev\abs{\xi_{\err+1}(t)}^p 
\le e^{Ct}. 
\label{zetabd9}\ee
Granting this for the 
moment, we apply H\"older's inequality
 to the other  expectations in \eqref{covbd3}, 
then \eqref{PABprob1} and the moment bounds \eqref{zetabd9}
to arrive at 
\be
\lvert \Cov(\om_{i_0}(t),\om_{\err+1}(t))\rvert \le 
C_1 e^{-C_2(\err-i_0)}. 
\label{covbd3.5}\ee

Now to verify \eqref{zetabd9} term by term.
For $\om_{i}(t)$  this is simply 
a finite moment under the $\mu^\te$ distribution.  
The other increment variables we express so that 
 Lemma \ref{hmombdlm} applies to each case. 
For the boundary increment $\om_{\err+1}(t)$
write 
\be
\begin{split}
\om_{\err+1}(t)&=h_{\err}(t)-h_{\err+1}(t)=
h_{\err}(t)-h_{\err}(0)+h_{\err}(0) -h_{\err+1}(0)\\
&= h_{\err}(t)-h_{\err}(0)+\om_{\err}(0)
\end{split}\label{omerraux}\ee
utilizing the frozen column $h_{\err+1}$ as in \eqref{omellincr}.
 Utilizing  \eqref{initzetaxi}
and \eqref{coupling3}
\begin{align*}
{\ze_{i}(t)}={h^\ze_{i-1}(t)-h^\ze_{i}(t)}
\le {h_{i-1}(t)}-{h_{i}(0)}
={h_{i-1}(t)}-{h_{i-1}(0)}+\om_i(0)
\end{align*} 
with a similar upper bound for $-\ze_i(t)$, which together give
\[
\abs{\ze_{i}(t)} \le 
\bigl({h_{i-1}(t)}-{h_{i-1}(0)}\bigr)
\vee\bigl({h_{i}(t)}-{h_{i}(0)}\bigr) +\abs{\om_i(0)}.
\]
For $\xi_{\err+1}(t)$ use stochastic monotonicity of columns
to begin with 
\[
{\xi_{\err+1}(t)}={h^\ze_{\err}(t)-h^\ze_{\err+1}(t)}
\le {h_{\err}(t)}-{h_{\err+1}(0)}
\]
and continue as in \eqref{omerraux}. 
The completes the verification of \eqref{zetabd9}. 

We have proved the  exponential decay of the covariance
in \eqref{covbd3.5}. 
The same arguments work for all the covariances on lines
\eqref{varhaux10a}--\eqref{varhaux10d}, and we state the lemma
in this generality.  The bounds extend also to the infinite-volume
stationary process because Lemma \eqref{hmombdlm} gives moment
bounds that ensure uniform integrability.   

\begin{lm} There exist constants $C_i=C_i(t,\te)\in(0,\infty)$ such
that, for all stationary $\lrt$ processes
and all $i,j\in\Zb$ and $s\in[0,t]$, 
\be
\Cov(\om_{i}(s),\om_{j}(t))  \le 
C_1 e^{-C_2\abs{i-j}}.\label{covbd4}\ee
The bound is also valid for the infinite-volume stationary 
process. 
\end{lm}  

Now we can complete the proof of the covariance formula.

\begin{proof}[Proof of Theorem \ref{varh1thm}]
  Since the bound \eqref{covbd4} hold uniformly 
as $-\ell,\err\to\infty$, the sums on lines 
\eqref{varhaux10b}--\eqref{varhaux10d} vanish
while the sums on line \eqref{varhaux10a} converge 
to
\[
\sum_{i\le0\,,\,j>z} 
\Cov(\om_i(0),\om_j(t)) 
\;+\; 
\sum_{i\le z\,,\,j>0} 
\Cov(\om_i(t),\om_j(0)).
\]
Translation invariance of  the stationary  infinite volume process
turns the above sum into the right-hand side of \eqref{varh1}. 
The left-hand side of \eqref{varhaux10a} converges to the left-hand
side of \eqref{varh1} as $-\ell,\err\to\infty$
 by Lemma \ref{lrtconvlm} and by the uniform
integrability given by Lemma \ref{hmombdlm}. 
\end{proof}

\bibliography{refsmarton}

\begin{thebibliography}{1}

\bibitem{valak}
M.~Bal{\'a}zs.
\newblock Microscopic shape of shocks in a domain growth model.
\newblock {\em J. Stat. Phys.}, 105(3/4):511--524, 2001.

\bibitem{fluct}
M.~Bal{\'a}zs.
\newblock Growth fluctuations in a class of deposition models.
\newblock {\em Ann. Inst. H. Poincar\'e Probab. Statist.}, 39(4):639--685,
  2003.

\bibitem{sokvalak}
M.~Bal{\'a}zs.
\newblock Multiple shocks in bricklayers' model.
\newblock {\em J. Stat. Phys.}, 117:77--98, 2004.

\bibitem{rwshscp}
M.~Bal{\'a}zs, Gy. Farkas, P.~Kov{\'a}cs, and A.~R{\'a}kos.
\newblock Random walk of second class particles in product shock measures.
\newblock {\em J. Stat. Phys.}, 139(2):252--279, 2010.

\bibitem{unipq3}
M.~Bal{\'a}zs, J.~Komj{\'a}thy, and T.~Sepp{\"a}l{\"a}inen.
\newblock Microscopic concavity and fluctuation bounds in a class of deposition
  processes.
\newblock {\em Annales de l'Institut Henri Poincar\'e. Probabilit\'es et
  Statistiques}, 48(1):151--187, 2012.

\bibitem{exists}
M.~Bal{\'a}zs, F.~Rassoul-Agha, T.~Sepp{\"a}l{\"a}inen, and S.~Sethuraman.
\newblock Existence of the zero range process and a deposition model with
  superlinear growth rates.
\newblock {\em Ann. Probab.}, 35(4):1201--1249, 2007.

\bibitem{convex}
M.~Bal{\'a}zs and T.~Sepp{\"a}l{\"a}inen.
\newblock A convexity property of expectations under exponential weights.
\newblock {\em {\tt http://arxiv.org/abs/0707.4273}}, 2007.

\bibitem{varj2nd}
M.~Bal{\'a}zs and T.~Sepp{\"a}l{\"a}inen.
\newblock Exact connections between current fluctuations and the second class
  particle in a class of deposition models.
\newblock {\em J. Stat. Phys.}, 127(2):431--455, 2007.

\bibitem{norris}
J.~R. Norris.
\newblock {\em Markov Chains}.
\newblock Cambridge Series in Statistical and Probabilistic Mathematics.
  Cambridge University Press, 1997.

\end{thebibliography}
\bibliographystyle{plain}
\end{document}